\PassOptionsToPackage{unicode}{hyperref}
\PassOptionsToPackage{hyphens}{url}
\PassOptionsToPackage{dvipsnames,svgnames,x11names}{xcolor}
\documentclass[
  12pt]{article}

\usepackage{amsmath,amssymb}
\usepackage{makecell}
\usepackage{iftex}
\ifPDFTeX
  \usepackage[T1]{fontenc}
  \usepackage[utf8]{inputenc}
  \usepackage{textcomp} 
\else 
  \usepackage{unicode-math}
  \defaultfontfeatures{Scale=MatchLowercase}
  \defaultfontfeatures[\rmfamily]{Ligatures=TeX,Scale=1}
\fi
\usepackage{lmodern}
\ifPDFTeX\else
\fi
\IfFileExists{upquote.sty}{\usepackage{upquote}}{}
\IfFileExists{microtype.sty}{
  \usepackage[]{microtype}
  \UseMicrotypeSet[protrusion]{basicmath} 
}{}
\makeatletter
\@ifundefined{KOMAClassName}{
  \IfFileExists{parskip.sty}{%
    \usepackage{parskip}
  }{
    \setlength{\parindent}{0pt}
    \setlength{\parskip}{6pt plus 2pt minus 1pt}}
}{
  \KOMAoptions{parskip=half}}
\makeatother
\usepackage{xcolor}
\makeatletter
\ifx\paragraph\undefined\else
  \let\oldparagraph\paragraph
  \renewcommand{\paragraph}{
    \@ifstar
      \xxxParagraphStar
      \xxxParagraphNoStar
  }
  \newcommand{\xxxParagraphStar}[1]{\oldparagraph*{#1}\mbox{}}
  \newcommand{\xxxParagraphNoStar}[1]{\oldparagraph{#1}\mbox{}}
\fi
\ifx\subparagraph\undefined\else
  \let\oldsubparagraph\subparagraph
  \renewcommand{\subparagraph}{
    \@ifstar
      \xxxSubParagraphStar
      \xxxSubParagraphNoStar
  }
  \newcommand{\xxxSubParagraphStar}[1]{\oldsubparagraph*{#1}\mbox{}}
  \newcommand{\xxxSubParagraphNoStar}[1]{\oldsubparagraph{#1}\mbox{}}
\fi
\makeatother

\usepackage{longtable,booktabs,array}
\usepackage{calc} 
\usepackage{etoolbox}
\makeatletter
\patchcmd\longtable{\par}{\if@noskipsec\mbox{}\fi\par}{}{}
\makeatother
\IfFileExists{footnotehyper.sty}{\usepackage{footnotehyper}}{\usepackage{footnote}}
\makesavenoteenv{longtable}
\usepackage{graphicx}
\makeatletter
\def\maxwidth{\ifdim\Gin@nat@width>\linewidth\linewidth\else\Gin@nat@width\fi}
\def\maxheight{\ifdim\Gin@nat@height>\textheight\textheight\else\Gin@nat@height\fi}
\makeatother
\setkeys{Gin}{width=\maxwidth,height=\maxheight,keepaspectratio}
\makeatletter
\def\fps@figure{htbp}
\makeatother

\makeatletter
\@ifpackageloaded{caption}{}{\usepackage{caption}}
\AtBeginDocument{%
\ifdefined\contentsname
  \renewcommand*\contentsname{Table of contents}
\else
  \newcommand\contentsname{Table of contents}
\fi
\ifdefined\listfigurename
  \renewcommand*\listfigurename{List of Figures}
\else
  \newcommand\listfigurename{List of Figures}
\fi
\ifdefined\listtablename
  \renewcommand*\listtablename{List of Tables}
\else
  \newcommand\listtablename{List of Tables}
\fi
\ifdefined\figurename
  \renewcommand*\figurename{Figure}
\else
  \newcommand\figurename{Figure}
\fi
\ifdefined\tablename
  \renewcommand*\tablename{Table}
\else
  \newcommand\tablename{Table}
\fi
}
\@ifpackageloaded{float}{}{\usepackage{float}}
\floatstyle{ruled}
\@ifundefined{c@chapter}{\newfloat{codelisting}{h}{lop}}{\newfloat{codelisting}{h}{lop}[chapter]}
\floatname{codelisting}{Listing}

\makeatother
\makeatletter
\@ifpackageloaded{caption}{}{\usepackage{caption}}
\@ifpackageloaded{subcaption}{}{\usepackage{subcaption}}
\makeatother

\ifLuaTeX
  \usepackage{selnolig}  
\fi
\usepackage[]{natbib}
\usepackage{bookmark}

\IfFileExists{xurl.sty}{\usepackage{xurl}}{} 
\hypersetup{
  pdftitle={Title},
  pdfauthor={Author 1; Author 2},
  pdfkeywords={3 to 6 keywords, that do not appear in the title},
  colorlinks=true,
  linkcolor={blue},
  filecolor={Maroon},
  citecolor={Blue},
  urlcolor={Blue},
  pdfcreator={LaTeX via pandoc}}

\newcommand{\anon}{1}
\usepackage{booktabs}
\usepackage{enumitem}
\newtheorem{theorem}{Theorem}

\newtheorem{proposition}[theorem]{Proposition}
\newtheorem{definition}[theorem]{Definition}

\newtheorem{property}{Property}
\newtheorem{condition}{Condition}
\newtheorem{example}{Example}

\begin{document}

\def\spacingset#1{\renewcommand{\baselinestretch}%
{#1}\small\normalsize} \spacingset{1}


\if1\anon
{
  \title{\bf The Pivotal Information Criterion}
  \author{Sylvain Sardy\thanks{
    The authors gratefully acknowledge the Swiss National Science Foundation grant 200021E\_213166.}\hspace{.2cm}\\
    Department of Mathematics, University of Geneva\\
    and \\
    Maxime van Cutsem\\
    Department of Mathematics, University of Geneva\\
    and \\
    Sara van de Geer \\
    Department of Mathematics, ETHZ}
  \maketitle
} \fi

\if0\anon
{
  \bigskip
  \bigskip
  \bigskip
  \begin{center}
    {\LARGE\bf Title}
\end{center}
  \medskip
} \fi

\bigskip
\begin{abstract}
The Bayesian and Akaike information criteria aim at finding a good balance between under- and over-fitting. They are extensively used every day by practitioners. Yet we contend they suffer from at least two afflictions: their penalty parameters $\lambda=\log n$ and $\lambda=2$ are too small, leading to many false discoveries, and their inherent (best subset) discrete
optimization is infeasible in high dimension. We alleviate these issues with the pivotal information criterion: PIC is defined as a continuous optimization problem, and the PIC penalty parameter $\lambda$ is selected at the detection boundary (under pure noise). PIC's choice of $\lambda$ is an upper quantile of a statistic that we prove to be (asymptotically) pivotal, provided the loss function is appropriately transformed.
As a result, simulations show a phase transition in the probability of exact support recovery with PIC, a phenomenon studied with no noise in compressed sensing.
Applied to real data, for similar predictive performances, PIC selects the least complex model among state-of-the-art learners.
\end{abstract}

\noindent%
\vfill

\newpage
\spacingset{1.8}

\section{Introduction}
As computing power and data storage capacity continue to increase, along with advances in optimization algorithms, researchers and practitioners allow themselves to consider highly flexible models with growing numbers of parameters.
Such flexibility may come at the expense of interpretability, parsimony, and sometimes scientific insight, as predictive performance alone becomes the primary objective.

We focus on model frugality and interpretability, using an information criterion to search for a parsimonious model and to control overfitting without relying on a validation set. To illustrate our objective, we consider the linear regression model
\begin{equation} \label{eq:linearmodel}
{\bf y}=\beta_0{\bf 1}+X{\boldsymbol \beta}+{\boldsymbol \epsilon},
\end{equation}
where $X$ is an $n\times p$ matrix, ${\boldsymbol \beta}\in\mathbb R^p$ is the parameter vector of interest, ${\boldsymbol \epsilon}$ are errors with variance  $\sigma^2$, and the intercept $\beta_0$ and $\sigma^2$ are nuisance parameters. A well-studied setting is the high-dimensional canonical linear model when $y_i=\beta_i+\epsilon_i$ for $i=1,\ldots,n$, with $n=p$ (i.e., $X=I_n$) and Gaussian errors. This simplified model has provided important insights into parameter estimation. In particular, the James--Stein estimator \citep{Jame:Stei:esti:1961}, a regularized maximum likelihood estimator, demonstrated that the maximum likelihood estimator $\hat {\boldsymbol \beta}^{\rm MLE}={\bf y}$ is inadmissible when $p> 2$.

As the number of predictors grows, it has become natural to assume that ${\boldsymbol \beta}$ is sparse, that is, that only a small fraction of the inputs carry predictive information. The parameter vector~${\boldsymbol \beta}$ is called sparse if
\begin{equation} \label{eq:needles}
{\cal S}=\left\{j\in\{1,\ldots,p\}: \beta_j \neq 0\right\}
\end{equation}
has a small cardinality $s:=|{\cal S}|$. For a linear model~\eqref{eq:linearmodel}, ${\cal S}$ identifies the inputs with predictive information, often referred to as ``needles in a haystack''. The resulting paradigm centers on detecting these needles, a task closely related to testing.

\subsection{Detection boundary and phase transition}

Within this sparse paradigm, \citet{Dono94b} in regression, and  \citet{Ingster1997} and \citet{HCDonoho2004} in hypothesis testing, showed the existence of a \emph{detection boundary} below which signal recovery is hopeless. Specifically, if the nonzero entries of ${\boldsymbol \beta}$ are below $\sigma \sqrt{2r\log n}$, where $r\in(0,1)$ is a sparsity index, then the nonzero entries are undetectable, otherwise detection becomes possible.
This observation suggests setting to zero the observations whose magnitude falls below $\sigma\sqrt{2r\log n}$, since attempting to estimate such components leads to excessive false discoveries. In particular, in wavelet smoothing, setting the threshold to $\sigma \sqrt{2 \log n}$ led to minimax properties \citep{Dono94b,Dono95asym}. \citet{IngsterX2010} extended the results to high-dimensional regression when the matrix $X$ is no longer the identity.

A phase transition analogous to the detection boundary in testing has also been observed in compressed sensing \citep{DonohoDL06,Candes:Romberg:2006}, corresponding to the noiseless linear model ${\bf y}=X{\boldsymbol \beta}$, where $X$ is $n\times p$ with $p>n$ and ${\boldsymbol \beta}$ is assumed $s$-sparse. For small sparsity levels and moderate dimensions, by using the $\ell_1$-minimization estimator solution to $\min_{{\boldsymbol \beta}\in{\mathbb R}^p} \| {\boldsymbol \beta} \|_1$ subject to ${\bf y}=X{\boldsymbol \beta}$, exact recovery occurs with probability close to one, but this probability sharply drops to zero as $s/p$ or $p/n$ increases beyond critical thresholds. LASSO \citep{Tibs:regr:1996} extends this approach to noisy settings by solving
\begin{equation}\label{eq:LASSO}
\min_{{\boldsymbol \beta}\in{\mathbb R}^p} \|{\bf y}-X{\boldsymbol \beta}\|_2^q+\lambda \| {\boldsymbol \beta} \|_1,
\end{equation}
 for $q=2$ and a regularization parameter $\lambda>0$ chosen to optimize the prediction risk via cross-validation, although this choice often fails to recover the support ${\cal S}$ \citep{earlyLASSO}.

\subsection{Information criteria} \label{subsct:IC}

The learner we consider here aims at a good estimation of ${\boldsymbol \beta}$, in particular at the detection of its nonzero entries, and is defined through an information criterion (IC). For a model parametrized by ${\boldsymbol \beta}\in{\mathbb R}^p$, and an intercept $\beta_0 \in {\mathbb R}$ and (possibly) a nuisance parameter $\sigma$, an IC takes the form
\begin{equation} \label{eq:IC0}
{\rm IC}=L({\boldsymbol \beta}, \beta_0, \sigma; {\cal D}) + \lambda C({\boldsymbol \beta}),
\end{equation}
where ${\cal D}=(X, \mathbf{y})\in\mathbb R^{n\times(p+1)}$ denotes the data, $L$ measures goodness of fit, and $C$ quantifies model complexity. The tuning parameter $\lambda>0$ is fixed a priori (i.e., without validation data or cross-validation) and governs the trade-off between $L$ and $C$ to select a model of appropriate dimension, thereby controlling under- and overfitting. Among information criteria, BIC \citep{Schw:esti:1978} is arguably the most widely used. It employs twice the negative log-likelihood for $L$ and the discrete complexity measure
\begin{equation} \label{eq:C0}
    C({\boldsymbol \beta})=\|\boldsymbol{\beta}\|_0=\#\{\beta_j\neq 0, \ j=1,\ldots,p \}
\end{equation}
which counts the number $s$ of active parameters. For BIC, the preset value is $\lambda=\log n$, whereas AIC adopts $\lambda=2$ \citep{AkaikeAIC73}. Minimizing either criterion yields estimators~$\hat {\boldsymbol \beta}$ in which some entries are exactly zero, facilitating interpretability by identifying influential inputs. The choice of $\lambda$ is crucial: smaller values tend to produce overly complex models that may overfit and include false positives, while larger values favor excessive sparsity.

A natural desideratum while minimizing the IC~\eqref{eq:IC0} is recovery of the truly relevant inputs of the linear model~\eqref{eq:linearmodel},
that is, matching $\hat {\cal S}=\{j\in\{1,\ldots,p\}: \hat \beta_j \neq 0\}$ with ${\cal S}$ in~\eqref{eq:needles}. This motivates studying the probability of exact support recovery,
\begin{equation} \label{eq:PESR}
{\rm PESR}={\mathbb P}(\hat {\cal S}={\cal S}),
\end{equation}
and its phase transition behavior as a function of $s=|{\cal S}|$, analogous to compressed sensing. This stringent property is not achieved by classical criteria such as BIC or AIC, even in ideal conditions such as independence of the inputs and high signal-to-noise ratio. Although dependence among columns of $X$ further complicates recovery, the independent-design setting provides an ideal benchmark for comparing sparse learners.

We use the Gaussian linear regression setting to illustrate the limitations of BIC, for which we propose a remedy in Appendix~\ref{app:BIC}.
Assuming Gaussian errors with variance $\sigma^2$,
\begin{equation*}
{\rm BIC}=\frac{1}{\sigma^2} {\rm RSS}_s+ \lambda s+ n \log \sigma^2 \quad {\rm with} \quad \lambda=\log n,
\end{equation*}
where the best size-$s$ residual sum of squares is
\begin{equation}\label{eq:RSSs}
{\rm RSS}_s:=\min_{\beta_0\in{\mathbb R}, {\boldsymbol \beta}\in\Omega_s }\|\mathbf{y}-\beta_0\mathbf{1}_n -X {\boldsymbol \beta}\|_2^2
\end{equation}
with $\Omega_s= \{ {\boldsymbol \beta} \in {\mathbb R}^p: \| {\boldsymbol \beta} \|_0=s\}$. Minimizing over $s\in\{ 0,1, \ldots, \min(n,p)\}$ yields a selected sparsity level $\hat s$. AIC replaces $\lambda=\log n$ with $\lambda=2$, leading to a less sparse model when $n>\exp(2)$. Regardless of whether one uses BIC, AIC, or some high-dimensional extensions such as \citet[EBIC]{ChenChen08} and \citet[GIC]{FanTang13}, these criteria depend on the nuisance parameter $\sigma^2$, whose estimation becomes challenging in high-dimensional settings ($p>n$); see, e.g., \citet{Fan:Guo:Hao:2012}.
\citet{BICGIC1984} proved that information criteria with $\lambda\to\infty$ and $\lambda/n\to 0$, hence BIC, are consistent. This result assumes a fixed number $p$ of candidate inputs as $n\to\infty$ and relies on a consistent estimate of $\sigma^2$ for overfitted models. These assumptions are not necessarily satisfied in our general setting, where $p$ may grow with $n$ and even exceed it. In particular, the BIC penalty $\log n$ does not account for the growing number of candidate models, motivating EBIC \citep{ChenChen08}.
Using the maximum likelihood estimate $\hat\sigma^{2,\rm MLE}={\rm RSS}_s/n$ yields the profile criterion (up to constants)
\begin{equation*}
{\rm BIC}=\log\frac{{\rm RSS}_s}{n} + s \frac{\log n}{n}.
\end{equation*}
Minimizing BIC is NP-hard in general, since computing ${\rm RSS}_s$ requires a search over all subsets of $s$ predictors. In practice, one relies on approximations such as greedy forward selection or the primal–dual active set (PDAS) method of \citet{2014InvPrL0}.

The discrete optimization behind BIC is certainly cumbersome, but, beyond the asymptotic regime of \citet{BICGIC1984}, the main reason why BIC fails to retrieve ${\cal S}$ in finite samples and moderate to high dimension stems from the choice $\lambda=\log n$. \citet{Schw:esti:1978} derived this term via a large-sample approximation to the Bayes factor, whose primary goal is model comparison rather than exact support recovery. Consequently, despite its widespread use, BIC is poorly suited for recovering ${\cal S}$, which in practice leads to suboptimal generalization and interpretability.
Section~\ref{sct:simu} illustrates this point empirically by comparing BIC to our proposed information criterion.

\subsection{Our proposal and paper organization}

We aim at constructing a learner that exhibits a phase transition analogous to compressed sensing, but in the presence of noise. Specifically, when $s/p$ or $p/n$ are small, the learner should recover ${\cal S}$ with high probability, while this probability should eventually drop to zero as these ratios increase and the nonzero components of ${\boldsymbol \beta}$ become undetectable.
To this end, we propose the \emph{pivotal information criterion} (PIC). Our main contributions are:
\begin{itemize}
\item \emph{A calibration principle for information criteria.} We pre-set $\lambda$ at the pivotal detection boundary, defined as the $(1-\alpha)$-quantile of the null law (under pure noise) of the zero-thresholding function,
 where the zero-thresholding function applied to data is defined as the smallest $\lambda$ for which the origin is a local minimizer. This requires extending the zero-thresholding function of \citet{CaroNickJairoMe2016} to information criteria with non-convex losses and non-convex penalties such as SCAD and MCP (Section~\ref{subsct:ztf}). The resulting learner retrieves ${\cal S}=\emptyset$ with probability $1-\alpha$ when $s=0$, whether or not the criterion is convex (Section~\ref{subsct:opti}).
\item \emph{A general recipe for pivotality.} The null law of the zero-thresholding function is made free of nuisance parameters by transforming the output ($\phi$) and the input ($g$) of the loss, or by using new loss functions. This new recipe contains square-root LASSO (which solves~\eqref{eq:LASSO} with $q=1$) as the special case of Gaussian errors with $\phi=\sqrt{\cdot}$  \citep{BCW11}, yields finite-sample pivotality for the location-scale family and asymptotic pivotality for the one-parameter exponential family (Section~\ref{subsct:phig}), and also applies to BIC (Appendix~\ref{app:BIC}).
\item \emph{A phase transition in the presence of noise.} The calibration at $s=0$ acts as a favorable anchor for the probability of exact support recovery, and we observe empirically a phase transition in this probability as $s$ increases. We do not prove existence of the phase transition here; in a companion paper, \citet{SaraPIC25} establish some theoretical results pointing to a phase transition for the location-scale family.
\item \emph{A plug-and-play learner.} With a continuous penalty, PIC requires neither cross-validation nor a combinatorial search such as forward selection, which becomes slow in high dimension: $\lambda$ is computed once from the design matrix, by simulation or in closed form (Sections~\ref{subsct:PDB} and~\ref{subsct:approx}), and a single fit returns the selected variables. On six real datasets, PIC selects the most parsimonious models among state-of-the-art learners, for a similar predictive performance (Section~\ref{subsct:real}).
\end{itemize}

The paper is organized as follows. Section~\ref{sct:PIC} defines PIC, extends the zero-thresholding function to non-convex criteria and derives the pivotal transformations $(\phi,g)$. Section~\ref{sct:computation} computes the pivotal detection boundary $\lambda_\alpha^{\rm PDB}$, assesses its sensitivity to the nuisance parameter for a finite sample, and computes the PIC learner. Section~\ref{sct:simu} illustrates the pivotality of PIC, studies its phase transition and its behavior with weak signals and correlated designs on synthetic data, and compares it with competitors on six real datasets. Section~\ref{sct:conclusion} concludes. Appendix~\ref{app:BIC} applies the PIC calibration to select $\lambda$ for BIC's discrete complexity measure, and the remaining appendices contain the proofs.

\section{The pivotal information criterion} \label{sct:PIC}

Many sparse learners such as Basis Pursuit \citep{CDS99}, LASSO, or Elastic Net \citep{ZouHastie05} rely on continuous penalties,  whereas AIC and BIC use a discrete complexity measure. However, what primarily distinguishes an information criterion \eqref{eq:IC0} from a sparse learner is that the regularization parameter $\lambda$ is \emph{pre-set}. Indeed, while learners select $\lambda$ from the data (e.g., via cross-validation  \citep{Arl_Cel:2010:surveyCV}, Stein unbiased risk estimation \citep{Stein:1981} or the slope heuristic \citep{BirgeMassart2007}), ICs rely on a closed-form model-based calibration (e.g., $\lambda=\log n$ for BIC or $\lambda=2$ for AIC). Two notable learners also use pre-set regularization parameters.
The first one is waveshrink \citep{Dono94b,Dono95asym}, which is LASSO~\eqref{eq:LASSO} in the canonical model $X=I_n$. The problem then decouples into $n$ univariate problems, each solved by soft-thresholding one observation $y_i$ with the universal threshold $\lambda=\sigma\sqrt{2\log n}$ setting all coefficients to zero with high probability under pure noise. In this univariate and sparse setting, the nuisance parameter $\sigma$ is easy to estimate. With a general design matrix $X$, and especially when $p>n$, estimating $\sigma$ is much harder.
The second notable learner is square-root LASSO \citep{BCW11}, which solves~\eqref{eq:LASSO} with $q=1$, that is, applies the transformation $\phi=\sqrt{\cdot}$ to LASSO's residual sum of squares. It is a particular instance of the quantile universal threshold of \citet{CaroNickJairoMe2016}, which selects $\lambda$ as the $(1-\alpha)$-quantile, under $H_0:{\boldsymbol \beta}={\bf 0}$, of the smallest $\lambda$ for which the origin ${\boldsymbol \beta}={\bf 0}$ is a minimizer of the criterion. For LASSO, the law of this statistic depends on $\sigma$. The square-root transformation removes this dependence and makes the statistic pivotal, so that its quantile can be computed without knowledge of $\sigma$. 

\subsection{Definition}\label{sct:PICdef}

The pivotal information criterion is a generalization of square-root LASSO beyond Gaussianity thanks to  the systematic construction of two transformations $\phi$ and $g$ tailored to the assumed probability law of the data. We establish here its general principle; a companion paper studies some theoretical properties  for the location-scale family and the $\ell_1$ penalty \citep{SaraPIC25}, and \citet{PIC4Cox2026} develop PIC for  survival analysis. Our key principle is to choose $\lambda$ at the detection boundary, formalized below.
\begin{definition}\label{def:PDB}
Let $\alpha\in(0,1)$ be a prescribed small probability. Under $H_0:{\boldsymbol \beta}={\bf 0}$, the $\alpha$-level detection boundary is
\begin{equation*}
    \lambda_\alpha=\inf\bigl\{\lambda>0:\ {\mathbb P}_{H_0}(\{{\boldsymbol \beta}={\bf 0} \mbox{ is a local minimizer of \eqref{eq:IC0} with } \lambda\})\geq 1-\alpha\bigr\},
\end{equation*}
the smallest penalty for which the origin is a local minimizer with probability at least $1-\alpha$ under pure noise. If $\lambda_\alpha$ does not depend on nuisance parameters, we call it the \emph{pivotal detection boundary} and denote it $\lambda_\alpha^{\rm PDB}$.
\end{definition}
We conjecture and observe empirically that choosing $\lambda<\lambda_\alpha$ leads to a high false discovery rate, whereas $\lambda>\lambda_\alpha$ unnecessarily compromises the chance of finding the needles. The universal threshold of Waveshrink satisfies such a property under $H_0$  \citep[equation (31)]{Dono94b} while being nearly minimax under $H_1$ (i.e., Besov spaces) \citep{Dono95asym}, for the orthogonal design.

Since discrete penalties lead to NP-hard problems, PIC instead relies on a class of continuous complexity measures $\mathcal{C}_{\ell_1}$  to facilitate optimization.
\begin{definition}\label{def:quasiL1}
Let $\rho:{\mathbb R}\to[0,\infty)$ be even, continuous, nondecreasing on ${\mathbb R}^+$ and $\rho(0)=0$. For $\boldsymbol{\beta}\in\mathbb R^p$, define $C(\boldsymbol{\beta})=\sum_{i=1}^p\rho(\beta_i)$. We say that $C\in\mathcal{C}_{\ell_1}$ is a \emph{first-order $\ell_1$-equivalent complexity penalty} if
$ \lim_{\epsilon\to0^+}\frac{C(\epsilon\boldsymbol{\beta})}{\epsilon}=\|\boldsymbol{\beta}\|_1$.
\end{definition}
This class of measures models complexity through the magnitude of the coefficients. It includes convex penalties (e.g., $\ell_1$) and nonconvex ones such as SCAD \citep{SCAD01} and MCP \citep{MCP2010} penalties. Although they may differ away from zero, all behave like the $\ell_1$ norm to first order near the origin, thereby promoting sparsity.
For BIC's nostalgics still wishing to use BIC's discrete complexity measure \eqref{eq:C0}, Appendix~\ref{app:BIC} shows how our method applies to find PIC's $\lambda$ for BIC.

The pivotal information criterion can now be defined.
\begin{definition}\label{def:PIC}
Given data $\cal D$ drawn from a parametric law, a base loss $l_n(\boldsymbol{\theta},\sigma;\mathcal{D})=\sum_{i=1}^n l(\theta_i,\sigma;\mathcal{D}_i)/n$, a complexity measure $C\in\mathcal{C}_{\ell_1}$, parameters $\boldsymbol{\beta}\in\mathbb R^p$, nuisance parameters $\beta_0\in\mathbb R$ and $\sigma$, and $\alpha\in(0,1)$, the \emph{Pivotal Information Criterion} is
\begin{equation}\label{eq:PIC}
\operatorname{PIC}=\phi\left(l_n({\boldsymbol \theta}, \sigma; \mathcal{D})\right) + \lambda_\alpha^{\rm PDB} C({\boldsymbol \beta}), \quad
{\boldsymbol \theta}=g(\beta_0{\bf 1}+X{\boldsymbol \beta}),
\end{equation}
where the pre-set penalty parameter $\lambda_\alpha^{\rm PDB}$ of Definition~\ref{def:PDB} is made pivotal thanks to the univariate functions $\phi$ and $g$ (applied componentwise), if they exist.
\end{definition}
The derivation of $\phi$ and $g$ is specific to the assumed law the data were generated from, and when $\phi$ and $g$ exist, their domain is also specific;  see Table~\ref{tab:PICsApplications}.
By construction, PIC takes the IC form defined in \eqref{eq:IC0} with composite loss $L=(\phi\circ l_n \circ g)$. Finding the functions $(\phi,g)$ is linked to the existence of a zero-thresholding function, which we discuss in Section~\ref{subsct:ztf}. Section~\ref{subsct:phig} derives $(\phi,g)$ for the location-scale and the one-parameter exponential families. The computation of $\lambda_\alpha^{\rm PDB}$ and of the PIC learner in practice is the topic of Section~\ref{sct:computation}.

\subsection{Zero-thresholding function} \label{subsct:ztf}

\citet{CaroNickJairoMe2016} derived the zero-thresholding functions of well-known sparse learners defined as solutions to convex optimization problems. The zero-thresholding function, when it exists, provides the smallest $\lambda$, denoted $\lambda_0$, for which the origin ${\boldsymbol \beta}=\mathbf{0}$ is a minimizer of the criterion.
For example, convex analysis leads to the zero-thresholding function  $\lambda_0=\| X^{\rm T} ({\bf y}- \bar y{\bf 1})\|_\infty$ for LASSO and $\lambda_0=\| X^{\rm T} ({\bf y}- \bar y{\bf 1})\|_\infty/\| {\bf y}- \bar y{\bf 1}\|_2$ for square-root LASSO. To observe the pivotality property of the latter, suppose that the data ${\bf y}$ are independent samples of a Gaussian random variable with mean $\mu$ and variance $\sigma^2$, then LASSO's zero-thresholding function is a random variable whose law depends on $\sigma$ (but not on $\mu$), as opposed to that of square-root LASSO that is independent of both $\mu$ and $\sigma$, and is therefore pivotal. LAD LASSO \citep{wang2007lad, Stest2025} and Rank LASSO \citep{Wang01102020} also have this pivotality property (along with being robust to non-Gaussian additive noise).

The pivotal information criterion~\eqref{eq:PIC} generalizes the pivotality of square-root LASSO and the existence of a zero-thresholding function to non-Gaussian models and non-convex learners. 
To that aim, we extend the definition of a zero-thresholding function to non-convex learners.
\begin{theorem}[Zero-thresholding function]\label{thm:ztf}
Consider the information criterion  \eqref{eq:IC0}  with $C\in {\cal C}_{\ell_1}$ and denote ${\boldsymbol \tau}=(\beta_0,\sigma)$ the nuisance parameters. Assume $L$ is twice differentiable in a neighbourhood of $(\mathbf 0,\hat{\boldsymbol \tau})$, where $\hat{\boldsymbol \tau}$ satisfies $\nabla_{{\boldsymbol \tau}} L (\mathbf 0,\hat {\boldsymbol \tau};\mathcal D)={\bf 0}$ and the Hessian $H_{\boldsymbol{\tau}}^L(\mathbf 0,\hat{\boldsymbol \tau}; \mathcal D)\succ0$. Then the zero-thresholding function
\begin{equation}\label{eq: ztf_eq}
\lambda_0(L,\mathcal D)=\left\|\nabla_{\boldsymbol{\beta}} L (\mathbf{0},\hat{\boldsymbol \tau};\mathcal D)\right\|_{\infty}
\end{equation}
satisfies:
\begin{enumerate}[label=(\roman*)]
    \item $( {\bf 0},\hat {\boldsymbol \tau})$ is a local minimizer of~\eqref{eq:IC0} if $\lambda>\lambda_0(L,\mathcal D)$;
    \item $( {\bf 0},\hat {\boldsymbol \tau})$ is not a local minimizer of~\eqref{eq:IC0} if $\lambda<\lambda_0(L,\mathcal D)$;
    \item If moreover $L$ is globally convex and $C$ is convex, then $( {\bf 0},\hat {\boldsymbol \tau})$ is a global minimizer of~\eqref{eq:IC0} if and only if $\lambda\geq \lambda_0(L,\mathcal D)$.
\end{enumerate}
\end{theorem}

Theorem~\ref{thm:ztf} characterizes the local minimality of the origin even for non-convex ICs~\eqref{eq:IC0}, whether the non-convexity stems from the penalty $C$ (e.g., SCAD or MCP) or from the loss $L$ (e.g., the Gaussian negative log-likelihood or its exponential; see Theorem~\ref{thm: loc-scale fam} below).
The $\alpha$-level detection boundary of Definition~\ref{def:PDB} can now be characterized. 
\begin{proposition}\label{prop:lamb under null}
 Let $\mathbf{Y}_0\sim p_{{\boldsymbol \tau}}$ denote observations generated under the pure-noise model ($\boldsymbol{\beta}=\mathbf{0}$ and ${\boldsymbol \tau}=(\beta_0,\sigma)$). Define the random variable
\begin{equation*}
    \Lambda = \lambda_0(L, (X, \mathbf{Y}_0))
\end{equation*}
and let $F_\Lambda$ denote its cumulative distribution function. Then $\lambda_\alpha = F_\Lambda^{-1}(1-\alpha):=\inf\{t: F_\Lambda(t)\geq 1-\alpha\}$. Moreover, with $\lambda=\lambda_\alpha$, the origin is a stationary point of~\eqref{eq:IC0} with probability $F_\Lambda(\lambda_\alpha)\geq 1-\alpha$, and, when $F_\Lambda$ is continuous,
$\mathbb{P}_{H_0}(\{{\boldsymbol{\beta}}={\bf 0}\text{ is a local minimizer of~\eqref{eq:IC0} with }\lambda=\lambda_\alpha\}) = 1 - \alpha$.
\end{proposition}
By Theorem~\ref{thm:ztf}(ii), the origin is a local minimizer only if $\Lambda\leq\lambda$, so the probability in Definition~\ref{def:PDB} is at most $F_\Lambda(\lambda)<1-\alpha$ for $\lambda<F_\Lambda^{-1}(1-\alpha)$. By Theorem~\ref{thm:ztf}(i), it is at least ${\mathbb P}(\Lambda<\lambda)\geq F_\Lambda(F_\Lambda^{-1}(1-\alpha))\geq 1-\alpha$ for $\lambda>F_\Lambda^{-1}(1-\alpha)$. Hence $\lambda_\alpha=F_\Lambda^{-1}(1-\alpha)$. When $F_\Lambda$ is continuous, ${\mathbb P}(\Lambda=\lambda_\alpha)=0$ and both bounds equal $1-\alpha$.

The value $\lambda_\alpha$ cannot be computed directly in practice because, being specific to the input matrix $X$, the distribution of the statistic~$\Lambda$  is unknown, and  may depend on the unknown nuisance parameters ${\boldsymbol \tau}=(\beta_0,\sigma)$ that generated the data at hand. Moreover, real datasets rarely arise from the pure-noise model, so that estimating $\beta_0$ or $\sigma$ reliably is challenging, especially in high dimension, which results in a poor calibration of $\lambda_\alpha$. By instead employing  PIC, whose transformed loss is $L(\boldsymbol\beta,\beta_0,\sigma;\mathcal D)=(\phi\circ l_n)\left(g(\beta_0\mathbf{1} + X\boldsymbol\beta),\sigma;\mathcal D\right)$, the law of $\Lambda$ can be made asymptotically pivotal through a zero-thresholding function involving $\phi'$ and $g'$. Here, we say that PIC is \emph{pivotal} if the law of $\Lambda$ does not depend on the nuisance parameters ${\boldsymbol\tau}$, and \emph{asymptotically pivotal} if there exist distribution functions $G_n$, depending on $X$ only, such that $\sup_{t\in\mathbb R}|\mathbb P_{\boldsymbol\tau}(\Lambda\le t)-G_n(t)|\to0$ for every fixed ${\boldsymbol\tau}$.

\begin{proposition}\label{prop: PIC-ztf}
Consider PIC \eqref{eq:PIC}, where $\phi$ and $g$ are strictly increasing differentiable univariate functions. Define the scalar-restricted loss $\tilde l_n(\theta,\sigma;\mathcal{D}):=l_n(\theta \mathbf 1_n,\sigma;\mathcal{D})$. Assume that under the null model ${\boldsymbol \beta}=\mathbf 0$ there exists $\hat{\boldsymbol\tau}=(\hat\beta_0,\hat\sigma)$ such that $g(\hat\beta_0)=\hat\theta^{\rm MLE}$ and $\hat\sigma=\hat\sigma^{\rm MLE}$, and that the scalar-restricted Hessian satisfies $\nabla^2_{(\theta,\sigma)}\tilde l_n(\hat\theta^{\rm MLE},\hat\sigma;\mathcal D)\succ 0$. Then PIC's zero-thresholding function  is
\begin{equation*}
\lambda_0\left(\phi\circ l_n \circ g, \mathcal{D}\right)
=\phi'\left(\tilde l_n\left(\hat\theta^{\rm MLE}, \hat \sigma; \mathcal{D}\right)\right) \cdot g'\bigl(g^{-1}(\hat\theta^{\rm MLE})\bigr) \cdot \left\|X^{\rm T}\nabla_{\boldsymbol{\theta}}l_n\left(\hat\theta^{\rm MLE}{\bf 1}, \hat\sigma; \mathcal{D}\right)\right\|_\infty.
\end{equation*}
\end{proposition}
Theorem~\ref{thm:ztf} requires $H_{\boldsymbol\tau}^L(\mathbf 0,\hat{\boldsymbol\tau}; \mathcal{D})\succ 0$ for the composite loss $L=\phi\circ l_n\circ g$. A second-order differentiation with respect to $(\beta_0,\sigma)$, together with the score equations eliminating first-order terms and the separable character of $l_n$, shows that this condition reduces to positive definiteness of the Hessian of $\tilde l_n$ at $(\hat\theta^{\rm MLE},\hat\sigma)$. The expression for $\lambda_0$ then follows from the chain rule. Notably, $\lambda_0$ does not depend on the specific penalty $C$, provided $C\in\mathcal C_{\ell_1}$.

\subsection{Finding \texorpdfstring{$(\phi,\, g)$}{(phi, g)}}\label{subsct:phig}
The transformations $(\phi,g)$ play a central role in rendering $\Lambda$  pivotal (Propositions~\ref{prop:lamb under null}-\ref{prop: PIC-ztf}), enabling calibration of $\lambda$ at the detection boundary without accurate estimation of nuisance parameters.
The function $g$ transforms the input of the loss $l_n$ and is reminiscent of the link function of generalized linear models \citep{NW72} and of variance stabilizing transformations \citep{Barlett1947}. The function $\phi$ rather transforms the output of the loss. It is important to observe that $\phi\circ l_n \circ g$ employed by PIC is just a composite loss, but by no means does it model the mechanism that generated the data through the functions $g$ and $\phi$. Take the Poisson law, for instance: the data may be generated from a linear model with the exponential canonical link, yet PIC's composite loss may employ a function $g\neq \exp$  to make the choice of $\lambda$ pivotal. Similarly, it is common practice to fit a linear model with the $\ell_2$ loss even when the data are Poisson.

To find a pivotal pair $(\phi,g)$, we consider loss functions $l_n$  of two classes of distributions: the location-scale and the one-parameter exponential families.
We let ${\cal D}=(X,\mathbf y)$ be the data, where $\mathbf y$ is a realization of $(Y_1,\dots,Y_n)$.

Employing the exponential transformation for $\phi$, the following theorem states the finite-sample pivotality of PIC's penalty for the location-scale family.

\begin{theorem}[Location-scale family]\label{thm: loc-scale fam}
Assume $f$ is a fixed baseline density such that, conditionally on $\mathbf x_i$,
\begin{equation*}
 f_{Y_i}(y_i\mid\mathbf x_i;\theta_i,\sigma)
=\frac{1}{\sigma}f\!\left(\frac{y_i-\theta_i}{\sigma}\right),
\quad f(\cdot)>0 \ \text{on}\ \mathbb R,\quad \sigma>0,
\end{equation*}
where $\theta_i=g(\beta_0+\mathbf x_i^{\rm T}\boldsymbol\beta)$. Let $l_n(\boldsymbol\theta,\sigma;{\cal D})$ be the average negative log-likelihood satisfying the assumptions of Proposition~\ref{prop: PIC-ztf}. Then PIC is finite-sample pivotal for $(\phi,g)=(\exp,{\rm id})$.
\end{theorem}

We refer to the resulting estimator as \emph{exponential LASSO}. The exponential transformation implies that exponential LASSO solves
\begin{equation*}
\min_{{\boldsymbol \beta}, \beta_0, \sigma}
\big(1/f_{\bf Y}({\bf y};{\boldsymbol \beta},\beta_0,\sigma)\big)^{1/n}
+\lambda\|\boldsymbol \beta\|_1.
\end{equation*}
For the Subbotin family with exponent $r$ (see below), profiling over $\sigma$ yields $\exp\big(l_n(\boldsymbol\theta,\hat\sigma;{\cal D})\big)\propto \big(\frac{1}{n}\sum_i|y_i-\theta_i|^r\big)^{1/r}$. This means that the $r$th root is the transformation that makes the $\ell_r$ loss pivotal, as stated in the following property.

\begin{property}
For the Subbotin family (including the Gaussian case for $r=2$) with density $f_{\bf Y}(y;\mu,\sigma)\propto \sigma^{-1}\exp(-|(y-\mu)/\sigma|^r)$,
the pivotal transformation for the loss  $l_n(\boldsymbol\theta,\sigma;{\cal D})=\frac{1}{n}\sum_i|y_i-\theta_i|^r$ is $(\phi,g)=(u^{1/r},{\rm id})$.
\end{property}
Location-scale distributions other than the Subbotin include the Huber, Gumbel and logistic laws. Asymptotic convergence properties of PIC for this family are studied in \citet{SaraPIC25}, and Section~\ref{sct:simu} illustrates the resulting phase-transition behavior empirically.

We now consider the one-parameter exponential family for which, conditionally on $\mathbf{x}_i$,
$\log f_{Y_i}(y_i \mid \mathbf{x}_i; \theta_i)
=
\theta_i T(y_i) - d(\theta_i) + b(y_i)$,
where the canonical parameter is linked to the linear predictor through $\theta_i = g(\beta_0 + \mathbf{x}_i^{\rm T}\boldsymbol{\beta})$, the canonical link corresponding to $g=\rm{id}$. We develop two approaches. The first one, described in Theorem \ref{thm: 1-param exp fam}, uses the classical negative log-likelihood as loss function $l_n$ and searches for a pivotal link function $g$. The second one, described in Theorem \ref{thm:weight score loss}, introduces a new class of losses yielding a weighted gradient allowing pivotality, yet keeping $g=\rm{id}$.

\begin{condition} \label{cond:1}
The vectors of variables $\mathbf x_1,\dots,\mathbf x_n\in\mathbb R^p$ in $X$  are fixed, bounded by $|x_{ij}|\leq B_n$, and standardized so that $\sum_{i=1}^n \mathbf{x}_i=\mathbf{0}$ and $\sum_{i=1}^n x_{ij}^2=n$ for $j=1,\dots,p$.
\end{condition}

\begin{condition} \label{cond:2}
The number of variables $p$ is allowed to tend to infinity such that $B_n^4 (\log(pn))^7 = O(n^{1-\nu})$ for some $\nu>0$.
\end{condition}

\begin{theorem}[One-parameter exponential family]\label{thm: 1-param exp fam}
Let $l_n$ be the average negative log-likelihood, and assume that, under $H_0:{\boldsymbol \beta}={\bf 0}$, the canonical parameter $\theta_0=g(\beta_0)$ does not depend on $n$ and lies in the interior of the natural parameter space. Under Conditions~\ref{cond:1} and \ref{cond:2}, PIC is asymptotically pivotal for
\begin{equation*}
      (\phi, g)=({\rm id},\, \tilde d^{-1})\quad\text{with}\quad \tilde d(\theta)=\int^\theta \sqrt{d''(s)}\,\mathrm{d}s.
\end{equation*}
More precisely, $\sup_{t\in\mathbb R}|\mathbb P(\sqrt n\,\Lambda\le t)-\mathbb P(\|Z\|_\infty\le t)|\to0$, where $Z\sim N(\mathbf 0,\hat\Sigma_X)$ with $\hat\Sigma_X=\frac{1}{n}\sum_{i=1}^n\mathbf x_i\mathbf x_i^{\rm T}$.
\end{theorem}

The proof of Theorem~\ref{thm: 1-param exp fam} (Appendix~\ref{app:expfam}) shows where pivotality comes from. Under $H_0$, the transformation $g$ makes the first two moments of the gradient at the origin pivotal: its mean is $\mathbf 0$ and its covariance is $\hat\Sigma_X/n$, whatever $\beta_0$ (up to a factor tending to one when $\beta_0$ is estimated). Its higher moments still depend on $\beta_0$, so $\Lambda$ is not pivotal in finite samples. Conditions~\ref{cond:1} and~\ref{cond:2} ensure that the gradient is approximately Gaussian, even when $p$ grows with $n$, and a Gaussian law depends only on its first two moments. Condition~\ref{cond:2} requires that no observation carries a non-negligible share of a column of $X$, which excludes, for instance, indicators of rarely observed categories. For the exponential law, pivotality even holds in finite samples since $\Lambda=\|X^{\rm T}({\bf y}-\bar y{\bf 1})\|_\infty/(n\bar y)$ (see Table~\ref{tab:PICsApplications}) is invariant under linear transformation of the data.

\begin{example}[Bernoulli]
In canonical form, $\log f(y\vert\theta)=y\theta-\log{(1+e^\theta)}$ with $\theta=\log\frac{\mu}{1-\mu}$ and $d''(\theta)=\frac{e^\theta}{(1+e^\theta)^2}$. So $\tilde d(\theta)=2\arctan(e^{\theta/2})$ (up to a constant) and
\begin{equation*}
    \theta_i = g(\eta_i)=2\log\left(\tan\frac{\eta_i+\omega}{2}\right),\quad \eta_i=\beta_0+\mathbf{x}_i^{\rm T}\boldsymbol{\beta},
\end{equation*}
where $\omega$ is an arbitrary centering constant. Correspondingly, the Bernoulli parameters are
\begin{equation*}
\mu_i=\frac{1}{1+\exp\!\big(-2\log(\tan\frac{\eta_i+\omega}{2})\big)}=\sin^2\left(\frac{\eta_i+\omega}{2}\right).
\end{equation*}
Choosing $\omega=\pi/2$ gives $\mu(0)=1/2$, leading to $\mu_i=(\sin(\eta_i)+1)/2$ for $\eta_i\in(-\pi/2,\pi/2)$, $i=1,\ldots,n$.
\end{example}

The pivotal reparametrization simplifies the variance, but  may yield link functions that are monotone only on a restricted interval. For Bernoulli,  the parameters need to satisfy the linear inequality constraints $\|\beta_0\mathbf{1}+X\boldsymbol{\beta}\|_\infty<\pi/2$. This induces constrained optimization in $(\beta_0,\boldsymbol{\beta})$, which is undesirable for computation. This motivates the second approach.

\begin{theorem}[Weighted score loss]\label{thm:weight score loss}
Under the same assumptions as Theorem \ref{thm: 1-param exp fam}, PIC is asymptotically pivotal, in the same sense, for $(\phi, g)=({\rm id},\, {\rm id})$ when considering the loss
\begin{equation*}
    l_n(\boldsymbol{\theta}, \mathcal{D})= \frac{1}{n}\sum_{i=1}^n \int^{\theta_i}\frac{d'(s) -T(y_i)}{\sqrt{d''(s)}}\mathrm{d}s. 
\end{equation*}
\end{theorem}
Under $H_0$, the fitted intercept coincides with the maximum likelihood estimate (Appendix~\ref{app:expfam}), so both approaches lead to the same statistic $\Lambda$ (Table~\ref{tab:PICsApplications}). If one is interested in using another link $g\neq {\rm id}$, the proof of Theorem~\ref{thm:weight score loss} derives the corresponding weighted loss.

\begin{example}[Bernoulli]
Basic integration leads to the loss $l_n(\boldsymbol{\theta}, \mathcal{D})=\frac{1}{n}\sum_{i=1}^n 2y_ie^{-\theta_i/2} + 2(1-y_i)e^{\theta_i/2}$ with $\theta_i=\beta_0+\mathbf{x}_i^{\rm T}\boldsymbol{\beta}$, or, with $\mu_i=1/(1+e^{-\theta_i})$,
\begin{equation*}
    l_n(\boldsymbol{\mu}, \mathcal{D})=\frac{1}{n}\sum_{i=1}^n 2y_i\sqrt{\frac{1-\mu_i}{\mu_i}} + 2(1-y_i)\sqrt{\frac{\mu_i}{1-\mu_i}}.
\end{equation*}
\end{example}
The discerning reader will notice that, with $\tilde y=2y-1$ to map $\{0,1\}$ into $\{-1,1\}$, the Bernoulli weighted score loss equals $2\exp(-\tilde y\theta/2)$, i.e., (twice) the AdaBoost exponential loss \citep{AdaBoost2000}.
\begin{example}[Poisson]
In canonical form, $\log f(y\vert\theta)=y\theta -\exp (\theta)$ with $\theta=\log(\mu)$ and $d'(\theta)=d''(\theta)=\exp(\theta)$. The corresponding weighted loss is $l_n(\boldsymbol{\theta}, \mathcal{D})=\frac{1}{n}\sum_{i=1}^n 2y_ie^{-\theta_i/2}+2e^{\theta_i/2}$ with $\theta_i=\beta_0+\mathbf{x}_i^{\rm T}\boldsymbol{\beta}$, or, with $\mu_i=\exp(\theta_i)$,
\begin{equation*}
    l_n(\boldsymbol{\mu}, \mathcal{D})=\frac{1}{n}\sum_{i=1}^n 2\frac{y_i}{\sqrt{\mu_i}} + 2\sqrt{\mu_i}.
\end{equation*}
\end{example}

Table~\ref{tab:PICsApplications} summarizes  instances of PIC transformation pairs with the corresponding gradients~$\nabla$ entering the zero-thresholding statistic $\Lambda=\|\nabla\|_\infty$. Within the Subbotin family, $r=2$ recovers the square-root LASSO, while $r=1$ yields the LAD LASSO. Notably, the Huber specification gives rise to a new estimator combining robustness, pivotality, and sparsity. For the one-parameter exponential family, the link function  and the loss are expressed in terms of the mean parameter $\mu$. We also include the Cox proportional hazards model. 
\begin{table}
\centering
\renewcommand{\arraystretch}{1.0}
\begin{tabular}{l c c c c c c }
\toprule
\textbf{Law}
& \textbf{Loss $l_n$}
& \textbf{$g(u)$}
& \textbf{${\rm Dom}(g)$}
& \textbf{$\phi(v)$}
& \textbf{$\nabla$}
& $c$ \\
\midrule

\multicolumn{7}{c}{\textit{Location-scale family}}\\

Gaussian
& NLL
& $u$
& ${\mathbb R}$
& $\exp{v}$
& $\sqrt{\frac{2\pi e}{n}}
\frac{X^{\rm T}\mathbf y_\bullet}{\|\mathbf y_\bullet\|_2}$
& $2\pi e$\\

Gaussian
& MSE
& $u$
& ${\mathbb R}$
& $\sqrt{v}$
& $\frac{X^{\rm T}\mathbf y_\bullet}
{\sqrt{n}\|\mathbf y_\bullet\|_2}$
& 1 \\

Gumbel
& NLL
& $u$
& ${\mathbb R}$
& $\exp{v}$
& $\frac{e^{1+\bar z}}{n}X^{\rm T}
\left(e^{\mathbf z}-\mathbf{1}\right)$
& $e^{2(\gamma+1)}$ \\

Huber ($k>0$)
& NLL
& $u$
& $\mathbb R$
& $\exp{v}$
& $\frac{L_0}{n\hat\sigma}
X^{\rm T}\boldsymbol\psi_k(\mathbf z)$
& $c_k$\\

Subbotin ($r>1$)
& M$r^{\rm th}$E
& $u$
& ${\mathbb R}$
& $v^{1/r}$
& $\frac{X^{\rm T}
\left(|\mathbf y_\bullet|^{r-1}\odot
\operatorname{sign}(\mathbf y_\bullet)\right)}
{n^{1/r}\|\mathbf y_\bullet\|_r^{r-1}}$
& $\frac{r^{2-2/r}\Gamma(2-1/r)}{\Gamma(1/r)}$\\

Laplace ($r=1$)
& MAE
& $u$
& ${\mathbb R}$
& $v$
& $\frac{X^{\rm T}\operatorname{sign}(\mathbf y_\bullet)}{n}$
& 1\\

\multicolumn{7}{c}{\textit{One-parameter exponential family}}\\

Bernoulli
& NLL
& $\frac{1+\sin u}{2}$
& $(-\frac{\pi}{2},\frac{\pi}{2})$
& $v$
& $\frac{X^{\rm T}\mathbf y_\bullet}
{n\sqrt{\bar y(1-\bar y)}}$
& 1\\

Bernoulli
& WSL
& $\frac{1}{1+e^{-u}}$
& $\mathbb R$
& $v$
& $\frac{X^{\rm T}\mathbf y_\bullet}
{n\sqrt{\bar y(1-\bar y)}}$
& 1\\

Poisson
& NLL
& $\frac{u^2}{4}$
& ${\mathbb R}^+$
& $v$
& $\frac{X^{\rm T}\mathbf y_\bullet}
{n\sqrt{\bar y}}$
& 1 \\

Poisson
& WSL
& $\exp{u}$
& $\mathbb R$
& $v$
& $\frac{X^{\rm T}\mathbf y_\bullet}
{n\sqrt{\bar y}}$
& 1 \\

Exponential
& NLL
& $\exp{u}$
& ${\mathbb R}$
& $v$
& $\frac{X^{\rm T}\mathbf y_\bullet}
{n\bar y}$
& 1 \\

\multicolumn{7}{c}{\textit{Survival analysis}}\\

Cox model
& NLPL
& $u$
& ${\mathbb R}$
& $\sqrt v$
& $(\S)$
& $\frac{1}{4\log n}$ \\

\bottomrule
\end{tabular}

\caption{
Examples of PIC's transformation pairs $(\phi,g)$.
NLL stands for negative log-likelihood, WSL for weighted score loss,
MSE for mean squared error and MAE for mean absolute error.
For the one-parameter exponential family, $g$ maps the mean parameter.
Throughout, $\mathbf y_\bullet=\mathbf y-\hat\mu\mathbf 1_n$
denotes the response centered at the null-model MLE $\hat\mu$,
and $\mathbf z=\mathbf y_\bullet/\hat\sigma$,
where $\hat\sigma$ is the null-model scale MLE.
For Huber,
$\psi_k(z)=\operatorname{sign}(z)\min\{|z|,k\}$,
$L_0=\hat\sigma C_k\exp\{n^{-1}\sum_i\rho_k(z_i)\}$
and $c_k=C_k^2\exp\{2\mathbb E[\rho_k(Z)]\}$,
where $C_k=\int_{\mathbb R}e^{-\rho_k(z)}\,dz$ is the normalizing constant and Z follows the standardized Huber distribution. $(\S)$ NLPL is Cox's negative log-partial likelihood
\citep{PIC4Cox2026}.
}
\label{tab:PICsApplications}
\end{table}

\section{Calibration and computation of the PIC learner}\label{sct:computation}

Section~\ref{sct:PIC} defined PIC through a penalty parameter $\lambda_\alpha^{\rm PDB}$ calibrated at the pivotal detection boundary, and characterized the local minimality of the origin through the zero-thresholding function (Theorem~\ref{thm:ztf}). Turning PIC into a learner requires two further steps: computing $\lambda_\alpha^{\rm PDB}$ (Sections~\ref{subsct:PDB} and~\ref{subsct:approx}), and computing an estimate $\hat{\boldsymbol\beta}$ with an algorithm that preserves this calibration, even when PIC is not convex (Section~\ref{subsct:opti}).

\subsection{Computing the pivotal detection boundary} \label{subsct:PDB}

The pivotal detection boundary $\lambda_\alpha^{\rm PDB}$ is the $(1-\alpha)$-quantile of the random variable $\Lambda$ of Proposition~\ref{prop:lamb under null}. For the pairs $(\phi,g)$ of Section~\ref{subsct:phig}, $\Lambda$ is pivotal for the location-scale family (and also for the exponential law) and asymptotically pivotal for the one-parameter exponential family. Note that the law of $\Lambda$ depends on the design matrix $X$ but not on the observed response ${\bf y}$, so $\lambda_\alpha^{\rm PDB}$ is computed once for a given $X$ matrix of inputs.

A natural approach to evaluate $\lambda_\alpha^{\rm PDB}$ is by Monte Carlo simulation. One generates $M$ response vectors $\mathbf{y}_0^{(m)}$, $m=1,\ldots,M$, under the null model $\mathbf{Y}_0\sim p_{{\boldsymbol \tau}_{\rm MC}}$ for a chosen nuisance parameter ${\boldsymbol \tau}_{\rm MC}$, computes $\lambda_0^{(m)}=\lambda_0(L,(X,\mathbf{y}_0^{(m)}))=\|\nabla^{(m)}\|_\infty$ with the closed-form gradients of Table~\ref{tab:PICsApplications}, and takes the empirical $(1-\alpha)$-quantile of $\lambda_0^{(1)},\ldots,\lambda_0^{(M)}$.
In the finite-sample pivotal case, Theorem~\ref{thm: loc-scale fam} guarantees pivotality for the given $n\times p$ matrix $X$, regardless of the chosen value of ${\boldsymbol \tau}_{\rm MC}$; consequently in that case, the only source of uncertainty in estimating $\lambda_\alpha^{\rm PDB}$ is due to the finiteness of $M$, and the larger $M$, the more precise the estimate.
For the asymptotically pivotal case, one must choose an arbitrary nuisance parameter ${\boldsymbol \tau}_{\rm MC}$ for the Monte Carlo simulation, whose influence on the estimate of $\lambda_\alpha^{\rm PDB}$ vanishes asymptotically by Theorem~\ref{thm: 1-param exp fam}.

\begin{figure}
    \centering
    \includegraphics[width=\linewidth]{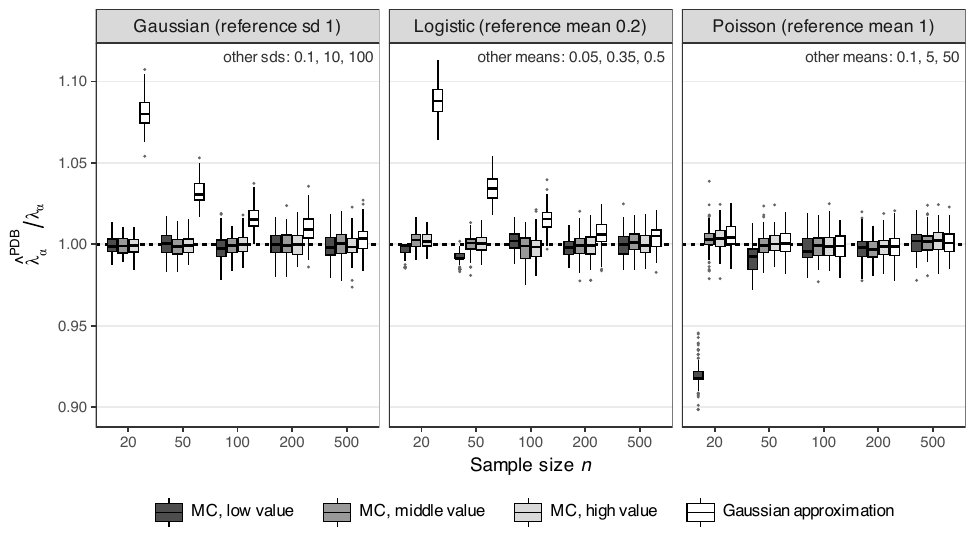}
    \caption{Sensitivity of the Monte Carlo calibration to the nuisance parameter, for the Gaussian, Bernoulli and Poisson models, with $M=2000$, $p=100$ and $\alpha=0.05$. Ratio of the calibrated $\hat\lambda_\alpha^{\rm PDB}$ to the reference calibration $\lambda_\alpha$ (dashed line), for Monte Carlo calibrations at three other values of the nuisance parameter (grey boxes) and for the Gaussian approximation of Proposition~\ref{prop:gauss} (white boxes), each repeated $100$ times.}
    \label{fig:calibration}
\end{figure}

To grasp the two sources of uncertainty decreasing when $M$ and $n$ get large,  Figure~\ref{fig:calibration} reports the results of the following empirical experiment. We consider  $\alpha=0.05$ and take a sequence of increasing sample sizes $n\in\{20,50,100,200,500\}$.
The design has $p=100$ independent standard Gaussian columns, standardized according to Condition~\ref{cond:1}; the data are either: (left) Gaussian with reference nuisance ${\boldsymbol \tau}:=(\mu,\sigma)=(0,1)$ (the benchmark, whose uncertainty is only due to $M$ since the Gaussian law belongs to the location-scale family); (middle) Bernoulli with reference nuisance $\tau:=\mu=0.2$; (right) Poisson with reference nuisance $\tau:=\mu=1$.
For these reference values $\tau$, we perform an extensive Monte Carlo simulation once with $10^5$ draws to calculate the three reference quantiles $\lambda_\alpha$.
Each is compared with calibrations at three other values (low, middle, high) of the nuisance parameter ($\mu=0, \sigma \in\{0.1,10,100\}$ for Gaussian, $\mu\in\{0.05,0.35,0.5\}$ for Bernoulli and $\mu\in\{0.1,5,50\}$ for Poisson data) using the default value of $M=2000$.
To evaluate the uncertainty due to a relatively small $M=2000$, we repeat the experiment $100$ times, calculating each time the ratio $\hat\lambda_\alpha^{\rm PDB}/\lambda_\alpha$, leading to the first three boxplots of Figure~\ref{fig:calibration} per sample size $n$.
We also add a fourth boxplot for the  Gaussian approximation of Proposition~\ref{prop:gauss} discussed below.

As expected, the finite-sample pivotality in the Gaussian case leads to boxplots aligned on the unit ratio of the estimated over the reference quantiles. The variability of the boxplots suggests that $M=2000$ Monte Carlo draws are sufficient. Similar conclusions apply to the asymptotically pivotal cases for large $n$, with a close alignment even for small $n$, except for an extreme value of the nuisance parameter with very few events (Poisson data with $\mu=0.1$ and $n=20$, about two events per simulated sample). The fourth boxplots show that the asymptotic Gaussian approximation can also be employed for moderate $n$, at the price of some conservativeness.
In conclusion, supported by the proven asymptotic pivotality property for the one-parameter exponential family and Figure~\ref{fig:calibration}, the estimation of PIC's $\lambda_\alpha^{\rm PDB}$ is not sensitive to the value of the nuisance parameter, as expected.
This motivates the following default values for our implementation (see Section~\ref{sct:reproduce}): $M=2000$ Monte Carlo draws, $\mu=1/2$ for the Bernoulli parameter and $\mu=5$ for the Poisson parameter.

\subsection{Approximating the pivotal detection boundary} \label{subsct:approx}

Despite its conceptual simplicity and highly parallelizable nature, Monte Carlo calibration may become computationally demanding, especially for large samples or when the evaluation of the gradient requires maximum likelihood computations without closed form (e.g., Gumbel). This motivates the following result, which bypasses explicit simulations of ${\bf Y}_0$.
\begin{proposition}[Asymptotic Gaussian calibration]\label{prop:gauss}
Under Conditions~\ref{cond:1} and~\ref{cond:2} (and, for the one-parameter exponential family, the assumptions of Theorem~\ref{thm: 1-param exp fam}), let $\hat \Sigma_X=\frac{1}{n}\sum_{i=1}^n \mathbf x_i\mathbf x_i^{\rm T}$, let $Z\sim N(0,c\,\hat\Sigma_X)$ with $c$
the constant of Table~\ref{tab:PICsApplications}, and let $q_n(\alpha)$ denote the $(1-\alpha)$-quantile of $\|Z\|_\infty/\sqrt n$.
For the location-scale and one-parameter exponential families of Table~\ref{tab:PICsApplications} with their transformations $(\phi,g)$,
\[
\rho_n:=\sup_{t\in\mathbb R}\Big|\,\mathbb P(\Lambda\le t)-\mathbb P\big(\|Z\|_\infty/\sqrt n\le t\big)\Big|\ \longrightarrow\ 0
\qquad (n\to\infty).
\]
Consequently, for every $\alpha\in(0,1)$,
\begin{itemize}
\item[(i)] (coverage) $\ \big|\mathbb P\{\Lambda\le q_n(\alpha)\}-(1-\alpha)\big|\le\rho_n\to0$;
\item[(ii)] (quantiles) $\ q_n(\alpha+\rho_n)\le F_\Lambda^{-1}(1-\alpha)\le q_n(\alpha-\rho_n)$ as soon as
$\rho_n<\min(\alpha,1-\alpha)$.
\end{itemize}
\end{proposition}

The convergence $\rho_n\to0$ is proven for the one-parameter exponential family in  Theorem~\ref{thm: 1-param exp fam} (Appendix~\ref{app:expfam}) and follows from \citet{SaraPIC25} for the location-scale family. Consequences~(i) and~(ii) then follow from the definition of $\rho_n$. In words, using the Gaussian quantile $q_n(\alpha)$ as the penalty level yields an asymptotically correct coverage, and the exact quantile of $\Lambda$ lies between two Gaussian quantiles whose levels differ from $\alpha$ by $\rho_n\to0$.

This asymptotic calibration reduces the problem to simulating Gaussian vectors with covariance $c\,\hat \Sigma_X$, thereby avoiding repeated maximum likelihood computations. A further simplification ignores the correlations among regressors: since the diagonal entries of $\hat\Sigma_X$ equal one under Condition~\ref{cond:1}, a union bound gives
\begin{equation*}
q_n(\alpha)\leq \sqrt{\frac{c}{n}}\,
\Phi^{-1}\!\left(1-\frac{\alpha}{2p}\right)
\leq \sqrt{
\frac{2c}{n}\log\!\left(\frac{2p}{\alpha}\right)
}.
\end{equation*}
This yields a simple closed-form rule for PIC's pre-set $\lambda$, which is conservative, and nearly exact when the regressors are weakly correlated.

\subsection{Computing the PIC learner} \label{subsct:opti}

With $\lambda=\lambda_\alpha^{\rm PDB}$, Proposition~\ref{prop:lamb under null} guarantees that, under $H_0$, the origin is a local minimizer of PIC with probability $1-\alpha$. This is a property of the criterion, not yet of an estimator. When PIC is strictly convex, local and global minima coincide, and Theorem~\ref{thm:ztf}(iii) characterizes exactly when the origin is the solution. When PIC is not convex, either because of the penalty (e.g., SCAD or MCP) or because of the transformed loss $\phi\circ l_n\circ g$, the criterion may have local minima. Nothing then guarantees that an algorithm finds the local minimum at the origin when it exists, so the calibration of $\lambda$ does not automatically transfer to the estimate. We therefore distinguish the local minima of PIC from the PIC learner $\hat{\boldsymbol\beta}$, which we define as the output of any solver satisfying three conditions that ensure this transfer.

\begin{definition}[Admissible solver and PIC learner]\label{def:PIClearner}
A solver for PIC~\eqref{eq:PIC} is \emph{admissible} if
(i) it is initialized at $({\bf 0},\hat{\boldsymbol\tau})$, where $\hat{\boldsymbol\tau}$ is the maximum likelihood estimate of the nuisance parameters under ${\boldsymbol\beta}={\bf 0}$;
(ii) it returns $({\bf 0},\hat{\boldsymbol\tau})$ whenever $({\bf 0},\hat{\boldsymbol\tau})$ is a stationary point of PIC;
(iii) its iterates converge to a stationary point of PIC.
The PIC learner $\hat{\boldsymbol\beta}$ is the ${\boldsymbol\beta}$-component of the output of an admissible solver.
\end{definition}
Stationarity is understood in the sense of the limiting subdifferential, which reduces to the usual first-order optimality condition for the $\ell_1$ penalty. Condition~(i) comes at no extra cost, since $\hat{\boldsymbol\tau}$ is needed anyway to evaluate the zero-thresholding function $\lambda_0(L,\mathcal D)$ (Proposition~\ref{prop: PIC-ztf}). Condition~(ii) requires the solver not to leave a stationary starting point; it holds for proximal gradient methods and their accelerated versions, as shown below for our algorithm. Condition~(iii) is the usual convergence requirement. The PIC learner is thus a stationary point of PIC, typically a local minimizer, but not necessarily its global minimizer. This choice is deliberate: the detection boundary of Definition~\ref{def:PDB} is defined through the local minimality of the origin, and any admissible solver returns the origin exactly when the origin is a local minimizer, as formalized below.

\begin{proposition}[Calibration of the PIC learner]\label{prop:PIClearner}
Let $C\in\mathcal{C}_{\ell_1}$ and assume the conditions of Theorem~\ref{thm:ztf}. Assume moreover that $\hat{\boldsymbol\tau}$ is the unique solution of $\nabla_{\boldsymbol\tau}L({\bf 0},{\boldsymbol\tau};\mathcal D)={\bf 0}$. Then, for any admissible solver and any data $\mathcal D$, the PIC learner satisfies
\begin{equation*}
\hat{\boldsymbol\beta}={\bf 0} \quad\Longleftrightarrow\quad \lambda_\alpha^{\rm PDB}\geq\lambda_0(L,\mathcal D).
\end{equation*}
Consequently, ${\mathbb P}_{H_0}(\hat{\boldsymbol\beta}={\bf 0})={\mathbb P}(\Lambda\leq\lambda_\alpha^{\rm PDB})\geq 1-\alpha$, with equality when $F_\Lambda$ is continuous.
\end{proposition}

Proposition~\ref{prop:PIClearner} follows from the stationarity condition of PIC at the origin. Since ${\boldsymbol\tau}$ is not penalized and the subdifferential of the penalty at ${\bf 0}$ is $[-1,1]^p$, a point $({\bf 0},{\boldsymbol\tau})$ is stationary if and only if ${\boldsymbol\tau}=\hat{\boldsymbol\tau}$ and $\|\nabla_{\boldsymbol\beta}L({\bf 0},\hat{\boldsymbol\tau};\mathcal D)\|_\infty\leq\lambda$, that is, $\lambda\geq\lambda_0(L,\mathcal D)$. Hence, if $\lambda_\alpha^{\rm PDB}\geq\lambda_0(L,\mathcal D)$, an admissible solver starts from a stationary point and returns it by~(ii); otherwise, no point with ${\boldsymbol\beta}={\bf 0}$ is stationary, and the solver cannot converge to one by~(iii). The probability statement follows since $\lambda_0(L,\mathcal D)$ has the law of $\Lambda$ under $H_0$.
This argument uses neither the convexity of PIC nor a specific solver: the calibration holds whether the non-convexity stems from the penalty or from the loss, and the solver only determines which nonzero stationary point is reached when $\lambda_\alpha^{\rm PDB}<\lambda_0(L,\mathcal D)$.

\paragraph*{Our algorithm.}
The PIC learner used for all experiments  is computed as follows.
\begin{enumerate}
\item \emph{Initialization.} Compute $\hat{\boldsymbol\tau}$ and $\lambda_0(L,\mathcal D)$.
\item \emph{$\ell_1$ path.} If $\lambda_\alpha^{\rm PDB}\geq\lambda_0(L,\mathcal D)$, set $K=1$ and $\lambda^{(1)}=\lambda_\alpha^{\rm PDB}$. Otherwise, consider a decreasing grid $\lambda_0(L,\mathcal D)=\lambda^{(1)}>\lambda^{(2)}>\cdots>\lambda^{(K)}=\lambda_\alpha^{\rm PDB}$. For $k=1,\ldots,K$, minimize PIC with the $\ell_1$ penalty and $\lambda=\lambda^{(k)}$ by FISTA with backtracking line search \citep{FISTA09}, initialized at $({\bf 0},\hat{\boldsymbol\tau})$ for $k=1$ and warm-started at the solution obtained for $\lambda^{(k-1)}$ for $k\geq2$.
\item \emph{Target-penalty stage.} If the target penalty $C$ is not the $\ell_1$ norm (e.g., SCAD or MCP), minimize PIC with penalty $C$ and $\lambda=\lambda_\alpha^{\rm PDB}$ by FISTA with backtracking line search, initialized at the output of step~2.
\end{enumerate}

This algorithm is admissible. Condition~(i) holds by construction. For condition~(ii), observe that a proximal gradient step maps a stationary point to itself (for non-convex $C$, provided the step size is small enough for the proximal map to be single-valued); the backtracking condition is then trivially satisfied, and the extrapolation term of FISTA vanishes since consecutive iterates coincide, so FISTA does not move from a stationary starting point. If $\lambda_\alpha^{\rm PDB}\geq\lambda_0(L,\mathcal D)$, the path reduces to $\lambda_\alpha^{\rm PDB}$, and $({\bf 0},\hat{\boldsymbol\tau})$ is stationary for the $\ell_1$ version of PIC as well, since the stationarity condition at the origin is the same for every $C\in\mathcal C_{\ell_1}$. Hence the $\ell_1$ stage returns the origin, and the target-penalty stage, initialized there, cannot leave it. Condition~(iii) is assumed for the last FISTA run; it is not guaranteed in general for accelerated methods applied to non-convex problems.

Neither the path nor the $\ell_1$ stage is needed for admissibility: FISTA applied directly to PIC from $({\bf 0},\hat{\boldsymbol\tau})$ is admissible as well. They matter when $\lambda_\alpha^{\rm PDB}<\lambda_0(L,\mathcal D)$. Following the $\ell_1$ solutions from the origin down to $\lambda_\alpha^{\rm PDB}$ with warm starts, as in \texttt{glmnet} \citep{GLMnet}, reduces the computational cost and reaches a nonzero stationary point in a stable way, and the $\ell_1$ solution provides the target-penalty stage with a good starting point.

\section{Simulation studies}\label{sct:simu}

\subsection{Illustration of (non-)pivotal transformations} \label{subsct:pivotalillustr}

\begin{figure}
    \centering
    \includegraphics[width=\linewidth]{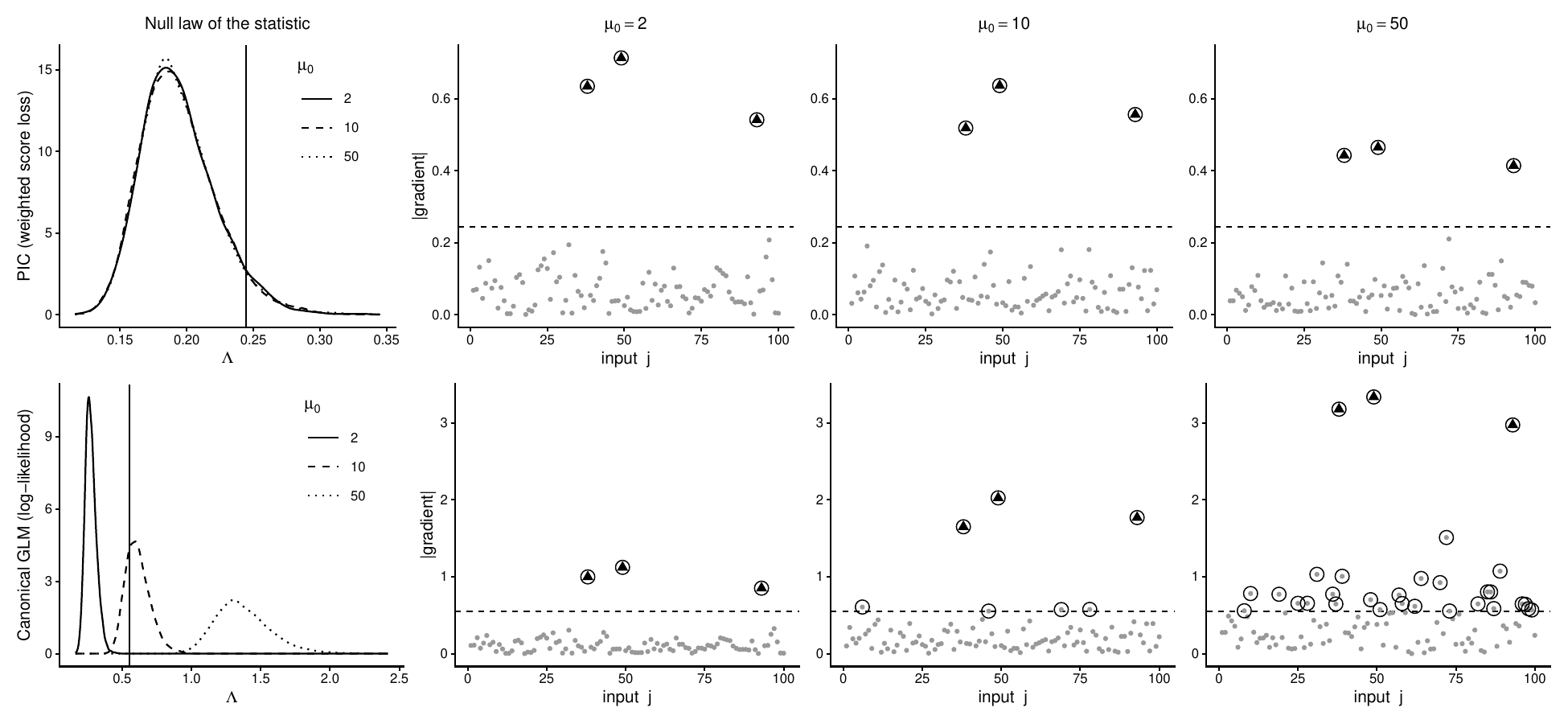}
    \caption{Pivotal versus non-pivotal calibration for Poisson data with $n=200$, $p=100$ orthogonal inputs and baseline mean $\mu_0\in\{2,10,50\}$. Top row: PIC (weighted score loss); bottom row: canonical GLM (negative log-likelihood, log link). Left: null law of the zero-thresholding statistic for each $\mu_0$, with $\lambda$ pre-set at the reference mean $5$ (vertical line). Other panels: magnitude of the gradient at the origin for one dataset with $s=3$ signals (triangles), with the same pre-set $\lambda$ (dashed line); circles mark the inputs above $\lambda$.}
    \label{fig:pivgrads}
\end{figure}

An information criterion pre-sets its penalty parameter, without knowing the nuisance parameter. Figure~\ref{fig:pivgrads} illustrates why this requires a pivotal statistic, by comparing, for Poisson data, the zero-thresholding statistic of PIC with that of the canonical generalized linear model. PIC uses the weighted score loss of Theorem~\ref{thm:weight score loss}, whose statistic is $\Lambda=\|X^{\rm T}({\bf y}-\bar y{\bf 1})\|_\infty/(n\sqrt{\bar y})$ (Table~\ref{tab:PICsApplications}); the canonical model uses the negative log-likelihood with the log link and $\phi={\rm id}$, whose statistic is $\|X^{\rm T}({\bf y}-\bar y{\bf 1})\|_\infty/n$. For both statistics, $\lambda$ is pre-set once, with the same procedure: the $(1-\alpha)$-quantile, $\alpha=0.05$, of the statistic simulated under $H_0$ with the default mean $\mu_{\rm MC}=5$ and $M=2000$ draws (see Section~\ref{subsct:PDB}).

The design has $n=200$ observations and $p=100$ orthogonal, centered and standardized columns ($X^{\rm T}X=nI_p$), so that the inputs are decoupled and the components of the gradient can be compared one by one with $\lambda$. We consider three baseline means $\mu_0\in\{2,10,50\}$, all different from $\mu_{\rm MC}$. The left panels display the null law of each statistic for the three values of $\mu_0$. The other panels display, for one dataset per value of $\mu_0$ with $s=3$ nonzero coefficients, the magnitude of each component of the gradient at the origin, $|\nabla_{\boldsymbol\beta}L({\bf 0},\hat\beta_0;\mathcal D)|$, together with the pre-set $\lambda$. The nonzero coefficients have random positions and signs, and magnitude $0.5/\sqrt{\mu_0}$, so that the signal-to-noise ratio is the same for the three means.

For PIC (top row), the three null laws coincide, as expected from pivotality, and the same pre-set $\lambda$ separates the three signals from the noise for every $\mu_0$. For the canonical model (bottom row), the null law is shifted and stretched by a factor $\sqrt{\mu_0}$, so that a single $\lambda$ cannot be valid for all values of the nuisance parameter. Calibrated at $\mu_{\rm MC}=5$, it is too large when $\mu_0<\mu_{\rm MC}$ and too small when $\mu_0>\mu_{\rm MC}$: the three signals are still detected, but $4$ and $27$ noise inputs exceed $\lambda$ for $\mu_0=10$ and $\mu_0=50$, respectively.

\subsection{Empirical phase transition analysis}  \label{subsct:empiricalpt}
We first study the probability of exact support recovery ($\mathrm{PESR}$) as a function of the sparsity level $s$, in settings where the only source of difficulty is the number of nonzero coefficients. To this end, the nonzero coefficients have a large, fixed magnitude and the inputs are independent, so that a drop in $\mathrm{PESR}$ can only be caused by a too large $s$. Increasing $n$ is then expected to push the transition towards larger values of $s$.

\paragraph*{Methods.}
We consider three versions of PIC: $\rm PIC{:}\ell_1$ and $\rm PIC{:}SCAD$, which use the $\ell_1$ and SCAD (with default $a=3.7$) penalties and are computed with the algorithm of Section~\ref{subsct:opti}, and $\rm PIC{:}\ell_0$, which uses BIC's discrete complexity measure~\eqref{eq:C0} with $\lambda$ calibrated at the pivotal detection boundary of the forward-selection zero-thresholding function (Appendix~\ref{app:BIC}); $\rm PIC{:}\ell_0$ is  run only for the Gaussian model. For PIC, $\alpha=0.05$ and $\lambda_\alpha^{\rm PDB}$ is computed once per design by Monte Carlo with $M=2000$ (Section~\ref{subsct:PDB}). For the Gaussian model, PIC uses the mean squared error loss with $\phi=\sqrt\cdot$ and $g={\rm id}$; for the Bernoulli and Poisson models, it uses the weighted score loss of Theorem~\ref{thm:weight score loss} with $g={\rm id}$ (Table~\ref{tab:PICsApplications}). We compare these approaches against three established baselines: BIC, EBIC, and GLMNet. The first two perform variable selection by minimizing
\begin{gather*}
    {\rm BIC}(s) = -2\log L_n(\hat\beta_0, \hat{\boldsymbol{\beta}}_s) + s\log n\quad {\rm and} \quad
    {\rm EBIC}_\gamma(s) =  {\rm BIC}(s) +2\gamma\log\binom{p}{s},
\end{gather*}
where $L_n$ denotes the likelihood function and $\gamma = 0.5$ (default). Both BIC and EBIC, as well as $\rm PIC{:}\ell_0$, are implemented via forward selection, capped at $40$ steps. The third baseline uses the \texttt{glmnet} package in \textsf{R} \citep{GLMnet} to compute the LASSO estimator with the regularization parameter chosen by 10-fold cross-validation (\texttt{lambda.min}).

\paragraph*{Design.}
We fix $p=100$. The design matrix $X$ has independent standard Gaussian entries, is standardized according to Condition~\ref{cond:1}, and is fixed for each model and sample size. For each sparsity level $s\in\{0,1,\ldots,25\}$, we generate $m=200$  datasets: the support of $\boldsymbol\beta$ is drawn uniformly random, and its $s$ nonzero entries have random signs and a fixed magnitude. The responses are generated from  $\eta_i=\beta_0+\mathbf x_i^{\rm T}\boldsymbol\beta$, $i=1,\ldots,n$, as follows:
\begin{itemize}
\item Gaussian model: $y_i=\eta_i+\epsilon_i$ with $\epsilon_i\sim N(0,1)$, $\beta_0=0$, $|\beta_j|=3$;
\item Bernoulli model: $y_i\sim{\rm Bernoulli}(1/(1+e^{-\eta_i}))$, $\beta_0=0$ (baseline mean $0.5$), $|\beta_j|=2$;
\item Poisson model: $y_i\sim{\rm Poisson}(e^{\eta_i})$, $\beta_0=\log 5$ (baseline mean $5$), $|\beta_j|=0.5$.
\end{itemize}
Besides $\mathrm{PESR}$, we report the true positive rate $\mathrm{TPR}=\mathbb E[|\hat{\cal S}\cap{\cal S}|/|{\cal S}|]$ (set to one when $s=0$) and the false discovery rate $\mathrm{FDR}=\mathbb E[|\hat{\cal S}\setminus{\cal S}|/\max(|\hat{\cal S}|,1)]$. At $s=0$, the FDR is the probability of selecting at least one variable, which PIC calibrates at $\alpha$ (equivalently, the FWER).

\begin{figure*}
    \centering
    \includegraphics[width=\linewidth]{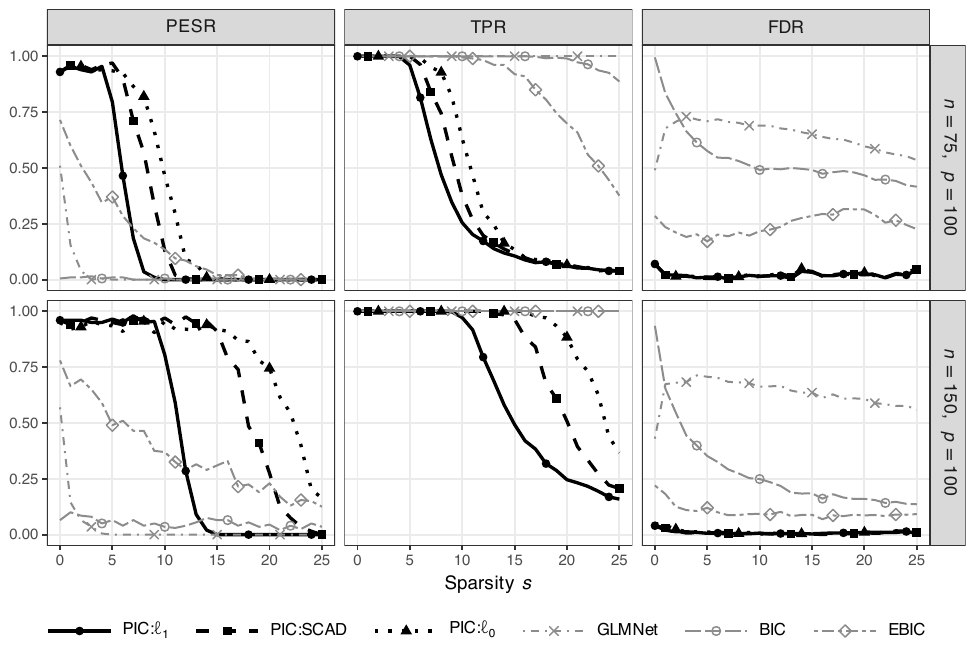}
    \caption{Gaussian phase transition analysis  with $p=100$ independent inputs and nonzero coefficients of magnitude $3$, for $n=75$ (top) and $n=150$ (bottom): $\mathrm{PESR}$ (left), TPR (middle) and FDR (right) as functions of the sparsity level $s$, averaged over $200$ datasets.}
    \label{fig:gaussian}
\end{figure*}

\paragraph*{Results.}
Figure~\ref{fig:gaussian} presents the results for the Gaussian model with $n\in\{75, 150\}$. Among the compared procedures, our three PIC-based methods are the only ones that exhibit a sharp phase transition in $\mathrm{PESR}$ as a function of $s$. In particular, the transition between near-perfect recovery and systematic failure is clearly visible, and shifts with increasing $n$. The non-convex penalties push the transition further than the $\ell_1$ penalty, whose shrinkage bias is known to hamper support recovery \citep{earlyLASSO}. The TPR and FDR panels show how the transition occurs: the FDR of PIC does not exceed $0.05$ for any $s$, and the drop in $\mathrm{PESR}$ is entirely due to missed detections. Beyond the transition, PIC thus becomes conservative rather than selecting noise variables. Among the baselines, EBIC consistently outperforms BIC, especially in the high-dimensional regime. GLMNet detects all true variables but at the price of many false ones (FDR between $0.54$ and $0.73$), so that its $\mathrm{PESR}$ is zero as soon as $s\ge1$.

\begin{figure*}
    \centering
    \includegraphics[width=\linewidth]{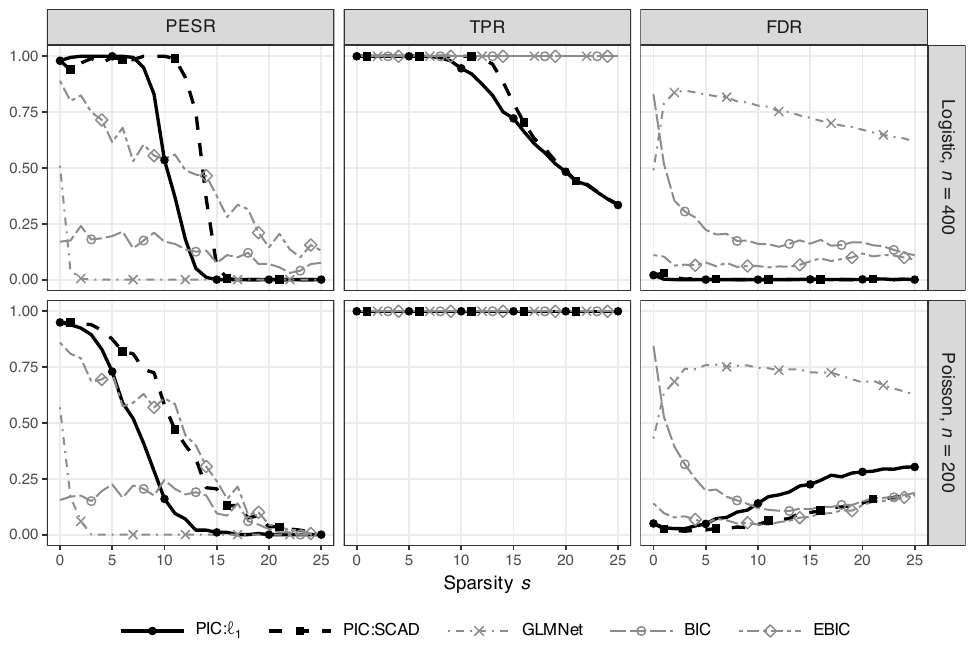}
    \caption{Bernoulli phase transition analysis  ($n=400$, nonzero coefficients of magnitude $2$, top) and the Poisson model ($n=200$, baseline mean $5$, nonzero coefficients of magnitude $0.5$, bottom), with $p=100$ independent inputs: $\mathrm{PESR}$ (left), TPR (middle) and FDR (right) as functions of the sparsity level $s$, averaged over $200$ datasets.}
    \label{fig:glm}
\end{figure*}

Figure~\ref{fig:glm} presents the results for the Bernoulli and Poisson models. For the Bernoulli model, the behavior mirrors that of the Gaussian model: our methods display a clear phase transition, whereas the baselines show a gradual degradation in performance. Compared to the Gaussian setting, the Bernoulli case highlights the increased difficulty of binary classification for exact support recovery and the need of a bigger dataset size $n=400$. For the Poisson model, the transition occurs earlier and is smoother, with $\mathrm{PESR}$ above $1/2$ up to $s=7$ for $\rm PIC{:}\ell_1$ and $s=10$ for $\rm PIC{:}SCAD$, comparable to EBIC ($s=11$). All methods detect all true variables in this setting, and the drop in $\mathrm{PESR}$ comes from false positives.

\subsection{Weak signals and correlated designs}  \label{subsct:weak}
The previous experiment isolates the effect of $s$. We now fix $s=5$ and make detection harder through signal amplitudes close to the detection boundary and through correlation among the inputs. The design has $p=100$ columns with Toeplitz correlation $\Sigma_{jk}=\rho^{|j-k|}$, $\rho\in\{0,0.5,0.8\}$, and is standardized and fixed for each model and value of $\rho$; the case $\rho=0$ isolates the effect of weak signals. The $s=5$ nonzero coefficients have random positions and signs and magnitudes $|\beta_j|=\sigma_{\rm eff}\,U_j$ with $U_j\sim U(0.5,1.5)$, where $\sigma_{\rm eff}$ is the inverse noise standard deviation at the baseline mean: $|\beta_j|\sim U(0.5,1.5)$ for the Gaussian model ($n=100$, $\epsilon_i\sim N(0,1)$, $\beta_0=0$), and $|\beta_j|\sim U(1,3)$ for the Bernoulli model ($n=300$, $\beta_0=0$, baseline mean $0.5$). Each design is also run with $s=0$ to check the calibration under correlation. The methods are those of Section~\ref{subsct:empiricalpt}, except $\rm PIC{:}\ell_0$, and $m=200$ datasets are generated per setting. Table~\ref{tab:weak} reports the TPR, the FDR and the mean model size $\bar{\hat s}$ at $s=5$, and, under $H_0$, the proportion of datasets in which at least one variable is selected,  whose target for PIC is $\alpha=0.05$.

Under the null model, PIC selects at least one variable in $4\%$ of the datasets for the Gaussian model and in $4\%$ to $9\%$ for the Bernoulli model, for every $\rho$, which is compatible with $\alpha=0.05$ given the Monte Carlo standard error (about $0.015$ for $200$ datasets). The competitors are far from this level: $10\%$ to $22\%$ for EBIC, $35\%$ to $43\%$ for GLMNet and $64\%$ to $96\%$ for BIC. With signals and $\rho\le0.5$, PIC selects models of size close to $s=5$ with an FDR between $0$ and $0.04$, at the price of missing some of the weakest signals (TPR between $0.88$ and $0.98$). The competitors detect almost all true variables, but with many false ones: GLMNet selects $19$ to $30$ variables (FDR $0.69$ to $0.82$), BIC has an FDR between $0.22$ and $0.46$, and EBIC, the closest competitor, between $0.07$ and $0.15$. With strong correlation ($\rho=0.8$), the FDR of $\rm PIC{:}\ell_1$ increases to $0.25$ for the Gaussian model and $0.16$ for the Bernoulli model, reflecting the known difficulty of the $\ell_1$ penalty with correlated designs, whereas $\rm PIC{:}SCAD$ is more robust (FDR $0.11$ and $0.05$) while keeping a TPR above $0.83$.

\begin{table}[htbp]
\centering
\resizebox{\linewidth}{!}{%
\begin{tabular}{lcccccccccccc}
\toprule
 & \multicolumn{4}{c}{$\rho = 0$} & \multicolumn{4}{c}{$\rho = 0.5$} & \multicolumn{4}{c}{$\rho = 0.8$} \\
\cmidrule(lr){2-5} \cmidrule(lr){6-9} \cmidrule(lr){10-13}
Method & TPR & FDR & $\bar{\hat s}$ & $H_0$ & TPR & FDR & $\bar{\hat s}$ & $H_0$ & TPR & FDR & $\bar{\hat s}$ & $H_0$ \\
\midrule
\multicolumn{13}{l}{\textit{Gaussian regression ($n = 100$)}} \\
PIC:$\ell_1$ & 0.91 & 0.01 & 4.6 & 0.04 & 0.88 & 0.04 & 4.7 & 0.04 & 0.81 & 0.25 & 5.6 & 0.04 \\
PIC:SCAD & 0.96 & 0.01 & 4.9 & 0.04 & 0.94 & 0.02 & 4.8 & 0.04 & 0.83 & 0.11 & 4.8 & 0.04 \\
GLMNet & 1.00 & 0.69 & 19.6 & 0.43 & 1.00 & 0.72 & 21.1 & 0.35 & 0.99 & 0.73 & 21.3 & 0.40 \\
BIC & 1.00 & 0.46 & 10.2 & 0.96 & 1.00 & 0.42 & 9.4 & 0.91 & 0.97 & 0.36 & 8.1 & 0.85 \\
EBIC & 1.00 & 0.13 & 5.9 & 0.22 & 1.00 & 0.14 & 6.0 & 0.17 & 0.96 & 0.15 & 5.8 & 0.17 \\
\midrule
\multicolumn{13}{l}{\textit{Bernoulli regression ($n = 300$)}} \\
PIC:$\ell_1$ & 0.90 & 0.00 & 4.5 & 0.07 & 0.89 & 0.01 & 4.5 & 0.04 & 0.82 & 0.16 & 5.1 & 0.09 \\
PIC:SCAD & 0.98 & 0.00 & 4.9 & 0.07 & 0.97 & 0.00 & 4.9 & 0.04 & 0.88 & 0.05 & 4.7 & 0.09 \\
GLMNet & 1.00 & 0.82 & 29.8 & 0.38 & 1.00 & 0.81 & 29.0 & 0.39 & 0.99 & 0.79 & 25.5 & 0.41 \\
BIC & 1.00 & 0.28 & 7.4 & 0.81 & 1.00 & 0.25 & 7.0 & 0.75 & 0.96 & 0.22 & 6.4 & 0.64 \\
EBIC & 1.00 & 0.08 & 5.5 & 0.14 & 1.00 & 0.07 & 5.4 & 0.10 & 0.95 & 0.10 & 5.3 & 0.12 \\
\bottomrule
\end{tabular}}
\caption{Weak signals and correlated designs ($s=5$, $p=100$, Toeplitz correlation $\rho^{|j-k|}$). For each $\rho$: TPR, FDR and mean model size $\bar{\hat s}$ at $s=5$, and at $s=0$ (under $H_0$) the proportion of datasets in which at least one variable is selected ($s=0$; target $\alpha=0.05$ for PIC). Averages over $200$ datasets.}
\label{tab:weak}
\end{table}

\subsection{Real data experiments} \label{subsct:real}
We now complement the synthetic analyses with an empirical study on real-world datasets, focusing on both variable selection and predictive performance. We consider six datasets: three regression tasks and three binary classification tasks. For each dataset, categorical variables are encoded using dummy variables when applicable, and missing values are handled via simple imputation using either the mean (for numerical variables) or the most frequent category (for categorical variables). The Diabetes data \citep{Efron2004LARS} include the $10$ baseline variables together with their $45$ two-way interactions and $9$ squares. Table~\ref{tab:real_datasets} summarizes the main characteristics of each dataset after preprocessing.

\begin{table}[ht]
\centering
\begin{tabular}{llccc}
\toprule
Task & Dataset & Source & $n$ & $p$ \\
\midrule
\textbf{Regression} & Prostate Cancer & Kaggle & 97 & 8 \\
&Diabetes & \textsf{R: lars} & 442 & 64 \\
&Riboflavin & \textsf{R: hdi} & 71 & 4088 \\
\midrule
\textbf{Binary classification} &Breast Cancer & UCI & 569 & 30 \\
&Ionosphere & UCI& 351 & 33 \\
& Sonar & UCI & 208 & 60 \\
\bottomrule
\end{tabular}
\caption{Considered datasets with some key characteristics.}
\label{tab:real_datasets}
\end{table}


For each dataset, we perform $100$ independent random train-test splits, where the test set represents $30\%$ of the observations and the remaining $70\%$ is used for training. On each split, the covariates are standardized with the mean and standard deviation of the training set, and all model fitting and variable selection procedures, including the computation of $\lambda_\alpha^{\rm PDB}$ ($\alpha=0.05$, $M=2000$), are carried out exclusively on the training data. The methods are those of Section~\ref{subsct:empiricalpt}, with $\rm PIC{:}\ell_0$ for the regression tasks only. Following variable selection, all methods are refitted without regularization on the selected covariates (Gaussian or Bernoulli regression), and predictive performance is evaluated on the test set. For classification, we report the test accuracy with a $0.5$ threshold on the predicted probability. For regression, we report the relative test mean squared error, that is, the test mean squared error divided by that of the intercept-only model fitted on the training set; it is unit-free and equals one for a model that does no better than the training mean.

Table~\ref{tab:real_datasets results} summarizes the results. On every dataset, a PIC procedure selects the smallest model. For regression, $\rm PIC{:}\ell_1$ achieves the lowest prediction error on Prostate Cancer and Riboflavin, and is within one standard error of the best method (BIC) on Diabetes, while selecting between $1.9$ and $6$ times fewer variables than GLMNet. On Riboflavin, where $p\gg n$, the advantage of PIC is clear: GLMNet selects $37.5$ variables on average and is unstable across splits, and BIC reaches the cap of $40$ forward steps. For classification, all accuracies  lie within $2$ percentage points. GLMNet is the most accurate on Ionosphere and Sonar, by $0.6$ and $1.5$ to $2$ points, but with about $4$ to $6$ times more variables than PIC, whereas $\rm PIC{:}\ell_1$ is the most accurate on Breast Cancer. Overall, for a similar predictive performance, PIC achieves greater parsimony, which makes it  appealing when interpretability matters.

\begin{table}[htbp]
\centering
\resizebox{\linewidth}{!}{%
\begin{tabular}{lcccccc}
\toprule
Dataset & PIC{:}SCAD & PIC{:}$\ell_1$ & PIC{:}$\ell_0$ & GLMNet & BIC & EBIC \\
\midrule
\multicolumn{7}{l}{\textbf{Regression} (relative test MSE)} \\
\midrule
\addlinespace
Prostate Cancer &
\makecell{$0.456 \pm 0.014$\\{\small $2.7 \pm 0.1$}} &
\makecell{$0.435 \pm 0.012$\\{\small $3.1 \pm 0.1$}} &
\makecell{$0.483 \pm 0.013$\\{\small $2.2 \pm 0.1$}} &
\makecell{$0.451 \pm 0.013$\\{\small $5.9 \pm 0.1$}} &
\makecell{$0.454 \pm 0.013$\\{\small $3.1 \pm 0.1$}} &
\makecell{$0.460 \pm 0.013$\\{\small $2.8 \pm 0.1$}} \\
\addlinespace
Diabetes &
\makecell{$0.535 \pm 0.005$\\{\small $4.7 \pm 0.1$}} &
\makecell{$0.531 \pm 0.005$\\{\small $4.6 \pm 0.1$}} &
\makecell{$0.548 \pm 0.006$\\{\small $3.5 \pm 0.1$}} &
\makecell{$0.536 \pm 0.007$\\{\small $17.6 \pm 0.7$}} &
\makecell{$0.529 \pm 0.006$\\{\small $6.3 \pm 0.1$}} &
\makecell{$0.542 \pm 0.005$\\{\small $4.7 \pm 0.1$}} \\
\addlinespace
Riboflavin &
\makecell{$0.411 \pm 0.013$\\{\small $6.0 \pm 0.2$}} &
\makecell{$0.411 \pm 0.013$\\{\small $6.0 \pm 0.2$}} &
\makecell{$0.602 \pm 0.023$\\{\small $2.6 \pm 0.1$}} &
\makecell{$1.323 \pm 0.425$\\{\small $37.5 \pm 1.1$}} &
\makecell{$0.754 \pm 0.039$\\{\small $40.0 \pm 0.0$}} &
\makecell{$0.492 \pm 0.022$\\{\small $8.0 \pm 1.0$}} \\
\addlinespace
\midrule
\multicolumn{7}{l}{\textbf{Classification} (test accuracy)} \\
\midrule
\addlinespace
Breast Cancer &
\makecell{$0.950 \pm 0.002$\\{\small $3.0 \pm 0.1$}} &
\makecell{$0.965 \pm 0.001$\\{\small $4.5 \pm 0.1$}} &
\makecell{--} &
\makecell{$0.963 \pm 0.001$\\{\small $12.9 \pm 0.2$}} &
\makecell{$0.962 \pm 0.001$\\{\small $5.2 \pm 0.1$}} &
\makecell{$0.963 \pm 0.001$\\{\small $4.3 \pm 0.1$}} \\
\addlinespace
Ionosphere &
\makecell{$0.867 \pm 0.003$\\{\small $4.3 \pm 0.1$}} &
\makecell{$0.872 \pm 0.003$\\{\small $4.6 \pm 0.1$}} &
\makecell{--} &
\makecell{$0.878 \pm 0.003$\\{\small $17.5 \pm 0.3$}} &
\makecell{$0.867 \pm 0.003$\\{\small $8.0 \pm 0.3$}} &
\makecell{$0.864 \pm 0.003$\\{\small $6.0 \pm 0.2$}} \\
\addlinespace
Sonar &
\makecell{$0.720 \pm 0.005$\\{\small $3.3 \pm 0.1$}} &
\makecell{$0.715 \pm 0.005$\\{\small $4.0 \pm 0.1$}} &
\makecell{--} &
\makecell{$0.735 \pm 0.004$\\{\small $19.8 \pm 0.6$}} &
\makecell{$0.729 \pm 0.005$\\{\small $7.6 \pm 0.3$}} &
\makecell{$0.722 \pm 0.005$\\{\small $4.4 \pm 0.2$}} \\
\addlinespace
\bottomrule
\end{tabular}}
\caption{Comparison of PIC with GLMNet, BIC and EBIC on regression and classification benchmarks. The first line of each entry reports the mean predictive performance: relative test mean squared error for regression (lower is better) and test accuracy for classification (higher is better). The second line reports the mean model size. Both are given as mean $\pm$ standard error over $100$ random train-test splits ($90$ for Riboflavin).}
\label{tab:real_datasets results}
\end{table}


\section{Conclusions and future work}\label{sct:conclusion}

BIC has been known to underpenalize complexity, causing false detections in many data analyses.
Attempts to alleviate this problem have been to increase BIC's penalty (e.g., EBIC and GIC) or to modify cross-validation with the one standard error rule \citep{cart84}.
These methods seem ad hoc and have not been designed to induce a phase transition. In empirical studies, these methods tend to select too many covariates.

The proposed method defines a new general framework to balance generalization and complexity in a way that mimics the phase transition achieved by compressed sensing.

The probabilistic selection of $\lambda$ is calibrated so that, under the complete null model, the probability of selecting at least one variable equals $\alpha$; this is a weak control of the family-wise error rate (which, under the complete null, coincides with the false discovery rate). Away from the null, the observed phase transition and control of the false discovery rate  require theoretical guarantees; \citet{SaraPIC25} initiated some results in that direction for the location-scale family.
PIC generalizes square-root LASSO to the location-scale and one-parameter exponential families, and can be naturally extended  beyond linear models and beyond regression or classification tasks;
 \citet{PIC4Cox2026} apply PIC to survival analysis.
The pivotality property of the gradient could also be employed further in gradient-based regularization.

\section{Reproducible research}\label{sct:reproduce}

The methods of this paper are implemented in the \textsf{R} package \texttt{picreg}, available on CRAN (\url{https://CRAN.R-project.org/package=picreg}); its source code, documentation and the \textsf{R} scripts reproducing all simulations are available at \url{https://github.com/VcMaxouuu/picreg}. The main function, \texttt{pic()}, fits the PIC learner for all the models 
 of Table~\ref{tab:PICsApplications}, with the $\ell_1$, SCAD or MCP penalty, using the algorithm of Section~\ref{subsct:opti}.

\section{Appendices}
\appendix

\section{Giving BIC a second chance} \label{app:BIC}

Section~\ref{subsct:IC} argues that
the Bayesian Information Criterion (BIC) has a low probability of exact support recovery in moderate to high dimension.
To fix BIC, we simply need to derive the zero-thresholding function of the discrete complexity measure~\eqref{eq:C0}  for an appropriate $(\phi,g)$-transformation such that the corresponding statistic $\Lambda$ is pivotal. As it turns out, the following theorem states that, with no transformation, BIC is already pivotal so that
 $\lambda$ can be selected at the pivotal detection boundary.

\begin{theorem}[Pivotal BIC]\label{thm: pivotal bic}
Consider the discrete measure of complexity $C$ defined in~\eqref{eq:C0}  and  let $l_n({\boldsymbol \theta},  \sigma; {\cal D})$ be twice the negative log-likelihood of the Gaussian law employed by BIC. Then, with the identity transformations $(g,\phi)=({\rm id},{\rm id})$, the zero-thresholding function is
\begin{equation}\label{eq: ztf BIC}
    \lambda_0(L,\mathcal{D})=\max_{s=1,\ldots,\min(p,n-2)} \frac{1}{s}\log\left(\frac{{\rm RSS}_0}{{\rm RSS}_s}\right)
\end{equation}
with ${\rm RSS}_s$ defined in~\eqref{eq:RSSs}, and the resulting IC is finite-sample pivotal. 
\end{theorem}

Theorem~\ref{thm: pivotal bic} shows that BIC can be embedded within the PIC framework
by calibrating its penalty at the pivotal detection boundary.
In practice, however, computing the BIC's zero-thresholding function \eqref{eq: ztf BIC} is infeasible in high dimension, as it requires solving an NP-hard optimization problem.
This drawback highlights another advantage of using PIC's continuous penalty.
A natural compromise for BIC is to replace best subset selection by a discrete notion of local minimality. For (one step) forward selection, the corresponding zero-thresholding function reduces to $\lambda_0(L,\mathcal{D} )=\log(\frac{{\rm RSS}_0}{{\rm RSS}_1})$. Since ${\rm RSS}_0/{\rm RSS}_1$ is pivotal under the null model, its pivotal detection boundary $\lambda_\alpha^{\rm PDB}$ is obtained exactly by the Monte Carlo calibration of Section~\ref{subsct:PDB}. This defines the $\rm PIC{:}\ell_0$ learner used in Section~\ref{sct:simu}.

\paragraph*{Proof of Theorem~\ref{thm: pivotal bic}.}
Twice the negative log-likelihood combined with the proposed transformations yields
\begin{equation*}
    L({\boldsymbol \beta}, \beta_0,\sigma;{\cal D})
=\frac{1}{\sigma^2}\frac{\|\mathbf{y}-\beta_0{\bf 1}_n-X{\boldsymbol \beta}\|_2^2}{n}
+\log(2\pi\sigma^2).
\end{equation*}
Let ${\rm RSS}_s$ be the best size-$s$ residual sum of squares defined in~\eqref{eq:RSSs}.
For fixed $s$, we have that $(\hat\sigma^{\rm MLE})^2={\rm RSS}_s/n$, so the profile information criterion is
${\rm IC}=1+\log(2\pi{\rm RSS}_s/n)+\lambda s$.
To ensure $s=0$ is a minimizer of the $\rm IC$, $\lambda$ must be large enough so that
\begin{equation*}
    1+\log(2\pi{{\rm RSS}}_0/n)\leq 1+\log(2\pi{{\rm RSS}}_s/n)+\lambda s, \quad \mbox{for all } s\geq 1.
\end{equation*}
Therefore the zero-thresholding function is that given in \eqref{eq: ztf BIC}.
So the statistic $\Lambda=\lambda_0(L,(X, {\bf Y}_0))$ is pivotal with respect to $(\beta_0,\sigma)$ since ${\rm RSS}_s/\sigma^2$ is pivotal for all $s$.

\section{Proof of Theorem~\ref{thm:ztf}}

Let $f({\boldsymbol \beta}, {\boldsymbol \tau})=L({\boldsymbol \beta}, {\boldsymbol \tau})+\lambda C({\boldsymbol \beta})$ with $C\in {\cal C}_{\ell_1}$.
Denote ${\boldsymbol\zeta}=(\boldsymbol{\beta}, {\boldsymbol \tau})$, ${\boldsymbol\zeta}^0=(\mathbf{0}, \hat{\boldsymbol \tau})$, where $ \hat{\boldsymbol \tau}$ minimizes $L$ at ${\boldsymbol \beta}={\bf 0}$.
Let ${\boldsymbol\zeta}^\epsilon={\boldsymbol\zeta}^0 +\epsilon {\bf h}$ with ${\bf h}=(d_{\boldsymbol\beta}, d_{{\boldsymbol \tau}})$. A second-order Taylor expansion of $L$ around ${\boldsymbol\zeta}^0$, combined with the first-order optimality conditions in ${\boldsymbol \tau}$, yields
\begin{equation*}
L({\boldsymbol\zeta}^\epsilon)=L({\boldsymbol\zeta}^0)+\epsilon d_{\boldsymbol{\beta}}^{\rm T}\nabla_{\boldsymbol{\beta}} L({\boldsymbol\zeta}^0)+\frac{\epsilon^2}{2}{\bf h}^{\rm T} H_{\boldsymbol\zeta}({\boldsymbol\zeta}^0) {\bf h}+ o(\epsilon^2).
\end{equation*}
Since $C\in\mathcal{C}_{\ell_1}$,
\begin{equation*}
\begin{aligned}
f({\boldsymbol\zeta}^\epsilon)-f({\boldsymbol\zeta}^0)
={}& \epsilon\, d_{\boldsymbol{\beta}}^{\rm T}
\nabla_{\boldsymbol\beta} L({\boldsymbol\zeta}^0)+ \frac{\epsilon^2}{2}\,
{\bf h}^{\rm T} H_{\boldsymbol\zeta}({\boldsymbol\zeta}^0){\bf h} + \lambda|\epsilon|\,\|d_{\boldsymbol{\beta}}\|_1
+ o(\epsilon^2).
\end{aligned}
\end{equation*}
By Hölder’s inequality, $|d_{\boldsymbol{\beta}}^{\rm T}\nabla_{\boldsymbol{\beta}} L({\boldsymbol\zeta}^0)|\leq \|\nabla_{\boldsymbol{\beta}} L({\boldsymbol\zeta}^0)\|_\infty\cdot\|d_{\boldsymbol{\beta}}\|_1$, which implies
\begin{equation*}
\begin{aligned}
f({\boldsymbol\zeta}^\epsilon)-f({\boldsymbol\zeta}^0)
\geq{}& |\epsilon|\left(\lambda
- \|\nabla_{\boldsymbol{\beta}} L({\boldsymbol\zeta}^0)\|_\infty\right)
\|d_{\boldsymbol{\beta}}\|_1 + \frac{\epsilon^2}{2}\,
{\bf h}^{\rm T} H_{\boldsymbol\zeta}({\boldsymbol\zeta}^0){\bf h}
+ o(\epsilon^2).
\end{aligned}
\end{equation*}
Writing ${\bf h}^{\rm T} H_{\boldsymbol\zeta}({\boldsymbol\zeta}^0){\bf h} = d_{\boldsymbol{\beta}}^{\rm T} H_{\boldsymbol{\beta}{\boldsymbol \beta}} d_{\boldsymbol{\beta}}+ 2 d_{\boldsymbol{\beta}}^{\rm T} H_{\boldsymbol{\beta}{\boldsymbol \tau}} d_\tau + d_{\boldsymbol \tau}^{\rm T} H_{{\boldsymbol \tau}{\boldsymbol \tau}} d_{\boldsymbol \tau}$, we distinguish directions with $d_{\boldsymbol{\beta}}=\mathbf 0$, for which positivity follows from $H_{{\boldsymbol \tau}{\boldsymbol \tau}}({\boldsymbol\zeta}^0)\succ 0$, and directions with $d_{\boldsymbol{\beta}}\neq\mathbf 0$, where the linear penalty term dominates all quadratic contributions for sufficiently small $\epsilon$.  Hence ${\boldsymbol\zeta}^0$ is a local minimizer of $f$ if
    $\lambda > \|\nabla_{\boldsymbol\beta} L({\boldsymbol\zeta}^0)\|_\infty$.
 On the contrary if $\lambda < \|\nabla_{\boldsymbol\beta} L({\boldsymbol\zeta}^0)\|_\infty$ then there is no local minimum at ${\boldsymbol\zeta}^0$. The case where the cost function in~\eqref{eq:IC0} is moreover globally convex was handled by \citet{CaroNickJairoMe2016} who proved a global minimum at ${\boldsymbol\zeta}^0$ if and only if  $\lambda \geq \|\nabla_{\boldsymbol\beta} L({\boldsymbol\zeta}^0)\|_\infty$.
This yields the stated zero-thresholding function.

\section{Proof of Theorem~\ref{thm: loc-scale fam}}
For the location-scale family,  the negative log-likelihood is $l_n(\theta,\sigma;  \mathcal{D})= \log \sigma - \frac{1}{n} \sum_{i=1}^n \log f((y_i-\theta)/\sigma)$. So with $\phi=\exp$, the zero-thresholding function is
\begin{align*}
\lambda_0\bigl(\exp(l_n), \mathcal{D}\bigr)
={}& \exp\!\left(
\log \hat \sigma^{\rm MLE}
- \frac{1}{n} \sum_{i=1}^n
\log f\!\left(
\frac{y_i-\hat \theta^{\rm MLE}}{\hat \sigma^{\rm MLE}}
\right)
\right) \cdot \left\|
\frac{1}{n\hat \sigma^{\rm MLE}}
\sum_{i=1}^n
\mathbf{x}_i\,\frac{
f'\!\left(
\frac{y_i-\hat \theta^{\rm MLE}}{\hat \sigma^{\rm MLE}}
\right)
}{
f\!\left(
\frac{y_i-\hat \theta^{\rm MLE}}{\hat \sigma^{\rm MLE}}
\right)
}
\right\|_\infty \\
={}& \frac{1}{n}\,
\exp\!\left(
- \frac{1}{n} \sum_{i=1}^n
\log f\!\left(
\frac{y_i-\hat \theta^{\rm MLE}}{\hat \sigma^{\rm MLE}}
\right)
\right) 
\cdot \left\|
\sum_{i=1}^n
\mathbf{x}_i\,\frac{
f'\!\left(
\frac{y_i-\hat \theta^{\rm MLE}}{\hat \sigma^{\rm MLE}}
\right)
}{
f\!\left(
\frac{y_i-\hat \theta^{\rm MLE}}{\hat \sigma^{\rm MLE}}
\right)
}
\right\|_\infty.
\end{align*}
Suppose ${\cal D}'=(X, {\bf y}')$ with ${\bf y}'=a{\bf y}+b$ with $a>0$, then, by definition of the location-scale family, ${\bf y}'$ is a sample from the same law as ${\bf y}$ but with parameters $(\theta', \sigma')$.
Moreover, by equivariance of the maximum likelihood estimator under affine transformations, we have that $\hat \theta'^{\rm MLE}=a\hat \theta^{\rm MLE}+b$ and $\hat \sigma'^{\rm MLE}=a\hat \sigma^{\rm MLE}$.
Hence the standardized residuals $(y_i'-\hat \theta'^{\rm MLE})/\hat \sigma'^{\rm MLE}=(y_i-\hat \theta^{\rm MLE})/\hat \sigma^{\rm MLE}$ are unchanged, so $\lambda_0\left(\exp \circ l_n, \mathcal{D}'\right)=\lambda_0\left(\exp \circ l_n, \mathcal{D}\right)$, proving that PIC is pivotal for the location-scale family.

\section{Proofs of Theorems~\ref{thm: 1-param exp fam} and~\ref{thm:weight score loss}}\label{app:expfam}
Assuming a one-parameter exponential family in canonical form $\log f(y_i\mid \theta_i)=\theta_i\, T(y_i)-d(\theta_i)+b(y_i)$, and a link $\theta_i=g(\eta_i)$ for $\eta_i=\beta_0+\mathbf{x}_i^{\rm T}\boldsymbol{\beta}$ leads to the negative log-likelihood
\begin{equation*}
    l_n(\boldsymbol{\theta}, \mathcal{D})=\frac{1}{n}\sum_{i=1}^n \bigl\{d(\theta_i) -\theta_i\,T(y_i)\bigr\}.
\end{equation*}
To create pivotality, we rather consider a general class of loss functions of the form
\begin{equation*}
    l_n(\boldsymbol{\theta}, \mathcal{D})=\frac{1}{n}\sum_{i=1}^n A(\theta_i) -w(\theta_i)\,T(y_i),
\end{equation*}
where $w$ is a differentiable function with $w'>0$. Imposing $A'(\theta)=d'(\theta)\,w'(\theta)$ ensures Fisher consistency, and choosing $w(\theta)=\theta$ recovers the negative log-likelihood. For $t\in\mathbb R$, let
\begin{equation*}
h_t(y,\mathbf x)=\bigl(T(y)-d'(g(t))\bigr)\,\kappa(t)\,\mathbf x,
\qquad \kappa(t)=w'(g(t))\,g'(t),
\end{equation*}
which is minus the gradient with respect to $\boldsymbol\beta$ of the loss of one observation $(y,\mathbf x)$ at $(\boldsymbol\beta,\beta_0)=(\mathbf 0,t)$, and let $P_nh_t=\frac{1}{n}\sum_{i=1}^n h_t(y_i,\mathbf x_i)$.

\emph{Null fit.} Under $H_0$, the derivative of the loss with respect to $\beta_0$ at $(\mathbf 0,t)$ is $\kappa(t)\{d'(g(t))-\bar T\}$, with $\bar T=\frac{1}{n}\sum_{i=1}^n T(y_i)$. In both theorems, $\kappa>0$ and $g$ is increasing. Hence, on the event $E_n$ that $\bar T$ lies in the interior of the mean space $d'(\Theta)$, $\Theta$ being the natural parameter space, the unique stationary point $\hat\beta_0$ satisfies $g(\hat\beta_0)=\hat\theta^{\rm MLE}=(d')^{-1}(\bar T)$ for any $w$; in particular, the intercept fitted with the weighted score loss of Theorem~\ref{thm:weight score loss} coincides with the maximum likelihood estimate. The second derivative at $\hat\beta_0$ equals $\kappa(\hat\beta_0)\,g'(\hat\beta_0)\,d''(\hat\theta^{\rm MLE})>0$, so Proposition~\ref{prop: PIC-ztf} applies and $\Lambda=\|P_nh_{\hat\beta_0}\|_\infty$ on $E_n$. Since $\bar T\to d'(\theta_0)$ almost surely, with $\theta_0=g(\beta_0)$ in the interior of $\Theta$, we have $\mathbb P(E_n^c)\to0$; on $E_n^c$, $\Lambda$ may be defined arbitrarily.

\emph{Factorization.} Thanks to the centering assumption $\sum_i\mathbf{x}_i=\mathbf{0}$, the term $d'(g(t))$ can be replaced by any constant, so that
\begin{equation*}
P_n h_{t}
=\kappa(t)\,\frac{1}{n}\sum_{i=1}^n
\bigl(T(y_i)-d'(g(t))\bigr)\mathbf{x}_i
=\kappa(t)\,\frac{1}{n}\sum_{i=1}^n
\bigl(T(y_i)-d'(g(\beta_0))\bigr)\mathbf{x}_i
\end{equation*}
for all $t$, in particular for $t=\hat\beta_0$ and for the true parameter value $t=\beta_0$. Therefore $P_n h_{\hat\beta_0}=r_n\,P_n h_{\beta_0}$ exactly, with $r_n=\kappa(\hat\beta_0)/\kappa(\beta_0)$. Define $\mathbf{v}_i=h_{\beta_0}(y_i,\mathbf{x}_i)=\xi_i\mathbf{x}_i$ with $\xi_i=\bigl(T(y_i)-d'(g(\beta_0))\bigr)\kappa(\beta_0)$, and $S_n=\frac{1}{\sqrt n}\sum_{i=1}^n\mathbf{v}_i$, so that $\sqrt n\,\Lambda=r_n\|S_n\|_\infty$ on $E_n$. The vectors $\mathbf{v}_i$ are independent, with $\mathbb E[\xi_i]=0$ and $\mathbb E[\xi_i^2]=\kappa(\beta_0)^2\,d''(g(\beta_0))$. Hence, imposing
\begin{equation*}
    w'(g(t))g'(t) = \frac{1}{\sqrt{d''(g(t))}}
\end{equation*}
gives $\mathbb E[\xi_i^2]=1$, so that $S_n$ has mean $\mathbf 0$ and covariance matrix $\hat\Sigma_X=\frac{1}{n}\sum_{i=1}^n \mathbf{x}_i\mathbf{x}_i^{\rm T}$, free of the nuisance parameter. Choosing $w(\theta)=\theta$ determines the pivotal link $g$ of Theorem~\ref{thm: 1-param exp fam}, whereas fixing the canonical link $\theta=g(\eta)=\eta$ leads to the loss of Theorem~\ref{thm:weight score loss}.
The first two moments of $S_n$ are thus exactly pivotal and, since a Gaussian law is determined by its first two moments, any Gaussian approximation of $S_n$ is pivotal. The remainder of the proof uses Conditions~\ref{cond:1} and~\ref{cond:2} to justify such an approximation for the maximum $\|S_n\|_\infty$ when $p$ may grow with $n$, and shows that replacing $\beta_0$ by $\hat\beta_0$, which only rescales $S_n$ by $r_n$, is asymptotically negligible.
In both cases, $\kappa(t)=1/\sqrt{d''(g(t))}$, so that $r_n=\{d''(\theta_0)/d''(\hat\theta^{\rm MLE})\}^{1/2}$. Since $d$ is infinitely differentiable on the interior of $\Theta$ with $d''>0$, the central limit theorem for $\bar T$ and the Delta method give $r_n=1+O_{\mathbb P}(n^{-1/2})$.

\emph{Gaussian approximation.} Since $\theta_0$ lies in the interior of $\Theta$, $T(Y)$ has a finite moment generating function in a neighborhood of zero, so the $\xi_i$ are independent and identically distributed sub-exponential random variables with mean zero and unit variance. Write $\|S_n\|_\infty=\max_{1\le k\le 2p}\tilde S_{nk}$, with $\tilde S_n=\frac{1}{\sqrt n}\sum_{i=1}^n\xi_i\tilde{\mathbf x}_i$ and $\tilde{\mathbf x}_i=(\mathbf x_i^{\rm T},-\mathbf x_i^{\rm T})^{\rm T}$. The coefficients $\tilde x_{ik}$ are non-stochastic, with $|\tilde x_{ik}|\le B_n$ and $\frac{1}{n}\sum_{i=1}^n\tilde x_{ik}^2=1$ by Condition~\ref{cond:1}, which also implies $B_n\ge1$. This setting is covered by the Gaussian approximation of \citet{Chernozhukov_2013}, whose rate condition is implied by Condition~\ref{cond:2} since $B_n\ge1$. Hence
\begin{equation*}
\rho_n':=\sup_{t\in\mathbb R}\bigl|\mathbb P(\|S_n\|_\infty\le t)-\mathbb P(\|Z\|_\infty\le t)\bigr|\to0,
\qquad Z\sim N(\mathbf 0,\hat\Sigma_X).
\end{equation*}

\emph{From $\beta_0$ to $\hat\beta_0$.} Let $\epsilon_n>0$ and $\delta_n=\mathbb P(|r_n-1|\,\|S_n\|_\infty>\epsilon_n)+\mathbb P(E_n^c)$. Since $\sqrt n\,\Lambda=r_n\|S_n\|_\infty$ and $\|S_n\|_\infty$ differ by at most $\epsilon_n$ outside an event of probability at most $\delta_n$,
\begin{equation*}
\mathbb P(\|S_n\|_\infty\le t-\epsilon_n)-\delta_n\le\mathbb P(\sqrt n\,\Lambda\le t)\le\mathbb P(\|S_n\|_\infty\le t+\epsilon_n)+\delta_n
\quad\text{for all } t\in\mathbb R.
\end{equation*}
The anti-concentration inequality of \citet{Chernozhukov_2015}, applied to the $2p$-dimensional Gaussian vector $(Z^{\rm T},-Z^{\rm T})^{\rm T}$, whose largest coordinate is $\|Z\|_\infty$ and whose coordinates have unit variance by Condition~\ref{cond:1}, gives $\sup_t\mathbb P(|\|Z\|_\infty-t|\le\epsilon)\le4\epsilon(\mathbb E\|Z\|_\infty+1)\le4\epsilon\{\sqrt{2\log(2p)}+1\}\le C\epsilon\sqrt{\log p}$ for a universal constant $C$ and $p\geq2$. Combined with the Gaussian approximation, this yields
\begin{equation*}
\sup_{t\in\mathbb R}\bigl|\mathbb P(\sqrt n\,\Lambda\le t)-\mathbb P(\|Z\|_\infty\le t)\bigr|\le\rho_n'+C\epsilon_n\sqrt{\log p}+\delta_n.
\end{equation*}
Moreover, $\|S_n\|_\infty=O_{\mathbb P}(\sqrt{\log p})$ since $\mathbb E\|Z\|_\infty\le\sqrt{2\log(2p)}$, and $r_n-1=O_{\mathbb P}(n^{-1/2})$. Choosing $\epsilon_n=\ell_n\sqrt{\log p/n}$ with $\ell_n\to\infty$ slowly enough gives $\delta_n\to0$, while $\epsilon_n\sqrt{\log p}=\ell_n\log p/\sqrt n\to0$ under Condition~\ref{cond:2}. Hence the Kolmogorov distance between $\sqrt n\,\Lambda$ and $\|Z\|_\infty$ tends to zero, where the law of $\|Z\|_\infty$ depends on $X$ only. This proves Theorems~\ref{thm: 1-param exp fam} and~\ref{thm:weight score loss}, as well as Proposition~\ref{prop:gauss} for the one-parameter exponential family, for which $c=1$.

\bibliography{article_bis}

\end{document}